# LONG-RANGE EXCLUSION PROCESSES, GENERATOR AND INVARIANT MEASURES

By Enrique D. Andjel and Hervé Guiol[1]

*Université de Provence and INP Grenoble*

We show that if $\mu$ is an invariant measure for the long range exclusion process putting no mass on the full configuration, $L$ is the formal generator of that process and $f$ is a cylinder function, then $Lf \in \mathbf{L}^1(d\mu)$ and $\int Lf \, d\mu = 0$. This result is then applied to determine (i) the set of invariant and translation-invariant measures of the long range exclusion process on $\mathbb{Z}^d$ when the underlying random walk is irreducible; (ii) the set of invariant measures of the long range exclusion process on $\mathbb{Z}$ when the underlying random walk is irreducible and either has zero mean or allows jumps only to the nearest-neighbors.

**1. Results and notation.** As the simple exclusion and zero-range processes, the long range exclusion process was introduced by Spitzer [8]. While the other two processes have been extensively studied (see [7] and the references therein), rather little is known about the long range exclusion process (see [3, 6, 9]). We believe that this is mostly due to the fact that the long range exclusion process does not satisfy the Feller property. The main goal of this paper is to, at least partially, overcome the difficulties related to the non-Feller character of the process.

The long range exclusion process, denoted by $(\eta_t)_{t \geq 0}$, is a continuous-time Markov process on the state space $\mathbf{X} = \{0,1\}^\mathbf{S}$, where $\mathbf{S}$ is a finite or countable set whose elements are called sites. If $\eta_t(x) = 1$, we say that at time $t$ there is a particle at $x$. Otherwise we say that $x$ is vacant at that time. Particles will attempt to move according to a Markov chain on $\mathbf{S}$ whose transition matrix is $p(x, y)$. The movement of particles will obey the following rules:

(i) *Exclusion rule*: There is always at most one particle at each $x \in \mathbf{S}$.

Received March 2003; revised November 2004.
[1]Supported in part by FAPESP Grant 96/04859-9.
*AMS 2000 subject classification.* 60K35.
*Key words and phrases.* Long range exclusion process, ergodic theorems, non-Feller process, invariant measures.

(ii) *Random clocks*: At each site there is a random clock which rings at times given by the jumps of a Poisson process of parameter 1. These processes are independent.

(iii) *Long jumps*: When a clock rings at an occupied site $x$, the particle at $x$ moves to site $X_\tau$, where $\{X_n\}_{n\geq 1}$ is a Markov chain on $\mathbf{S}$ with transition probability matrix $p$ and initial condition $X_0 \equiv x$, and $\tau$ is the first (positive integer) time $n$ such that $X_n$ is a vacant site (site $x$ itself is considered vacant during the jump; i.e., if the chain returns to site $x$ without visiting before an empty site, then the jump is cancelled). If the stopping time $\tau$ is infinite, then the particle disappears.

The construction of this process for an infinite number of particles is based on monotonicity arguments and is due to [6]. We recall his construction in the next section.

Let $\eta_t^\xi$ be the long range exclusion starting from the initial configuration $\xi$, let $S(t)$ be the semigroup, acting on the set of bounded measurable functions on $\mathbf{X}$, given by $S(t)f(\xi) = \mathbf{E}(f(\eta_t^\xi))$ and for any probability measure $\mu$ on $\mathbf{X}$ let $\mu S(t)$ be the unique probability measure such that

$$\int f \, d\mu S(t) = \int S(t) f \, d\mu,$$

for all bounded measurable $f$.

Denote by $\mathcal{I} = \{\mu : \mu S(t) = \mu \text{ forall } t > 0\}$ the set of invariant probability measures on $\mathbf{X}$ for the long range exclusion process. Obviously $\mathcal{I}$ is nonempty since it contains at least the measure concentrated on the empty configuration.

For any cylinder function $f$ and $\eta \in \mathbf{X}$ we write

$$
\begin{aligned}
Lf(\eta) = & \sum_{x,y \in \mathbf{S}} q(x,y,\eta) \eta(x)[1 - \eta(y)][f(\eta^{xy}) - f(\eta)] \\
& + \sum_{x \in \mathbf{S}} \delta(x,\eta)[f(\eta_x) - f(\eta)],
\end{aligned}
$$
(1)

where for $x \neq y$

$$q(x,y,\eta) = \mathbb{E}^x \left[ \prod_{i=1}^{\sigma(y)-1} \eta(X_i), \sigma(y) < \sigma(x), \sigma(y) < \infty \right],$$

$X_i$ is a Markov chain with transition matrix $p(x,y)$, the exponent $x$ on $\mathbb{E}$ denotes the starting point of the chain $X_i$,

$$\sigma(y) = \inf\{n \geq 1 : X_n = y\},$$

$$\delta(x,\eta) = \mathbb{E}^x \left[ \prod_{i=0}^{\infty} \eta(X_i), \sigma(x) = \infty \right],$$

$\eta^{xy}$ is configuration $\eta$ where the states of sites $x$ and $y$ have been interchanged:

$$\eta^{xy}(z) = \begin{cases} \eta(x), & \text{if } z = y, \\ \eta(y), & \text{if } z = x, \\ \eta(z), & \text{otherwise,} \end{cases}$$

and

$$\eta_x(z) = \begin{cases} \eta(z), & \text{if } z \neq x, \\ 0, & \text{if } z = x. \end{cases}$$

The operator $L$ can be thought of as the "formal" infinitesimal generator of the process. Note that the series above might not converge for some values of $\eta$ even if $f$ is a cylinder function. To see this, suppose the transition matrix $p(x,y)$ corresponds to a nearest-neighbor random walk on $\mathbb{Z}$ with drift toward the right. Then take an initial configuration $\eta \in \mathbf{X}$ such that $\eta(x) = 1 \ \forall \, x < 0$ and $\eta(0) = 0$ and apply the generator to $f(\eta) := \eta(0)$.

Our first result relates the operator $L$ to the invariant measures of our process: Denote by $\mathbf{1}$ the full configuration, that is, the element of $\mathbf{X}$ such that $\mathbf{1}(x) = 1$ for all $x \in \mathbf{S}$, and by $\nu_1$ the point mass measure on $\mathbf{1}$.

THEOREM 1.1. *Suppose that $p(\cdot, \cdot)$ is irreducible and $\nu_1 \in \mathcal{I}$. Let $\mu \in \mathcal{I}$ be such that $\mu(\{\mathbf{1}\}) = 0$ and for any finite subset $R \subset \mathbf{S}$ let $f_R(\eta) = \prod_{x \in R} \eta(x)$. Then, the series defining $Lf_R(\eta)$ converges $\mu$ a.e. Moreover, $Lf_R(\eta) \in \mathbf{L}^1(\mu)$ and*

$$\int L f_R \, d\mu = 0.$$

REMARK. Although we believe that the conclusions of this theorem also hold when $\nu_1 \notin \mathcal{I}$, this will require a different proof. Since we will apply the theorem in cases where $\nu_1 \in \mathcal{I}$ we have only treated this case.

In the last two sections of this paper we will use coupling techniques to determine the set of invariant measures for some long range exclusion processes on $\mathbb{Z}$. For this purpose we will need a process on $\mathbf{X} \times \mathbf{X}$ whose marginals are versions of the long range exclusion process and such that particles of different marginals move together as much as possible. The exact meaning of this will become clear in the next section. We denote by $\widetilde{L}$ the "formal" coupled generator. The expression of this generator, written below, is quite involved and we suggest the reader skip it until he has seen the next section:

$$\widetilde{L} f(\eta, \xi)$$

$$
\begin{align}
(2) \quad &= \sum_{\substack{\eta(x)=\xi(x)=1 \\ \eta(y)=\xi(y)=0}} q(x,y,\eta\xi)[f(\eta^{xy},\xi^{xy}) - f(\eta,\xi)] \\
(3) \quad &+ \sum_{\substack{\eta(x)=\xi(x)=1 \\ \eta(y)=\xi(z)=0}} q(x,y,z,\eta,\xi)[f(\eta^{xy},\xi^{xz}) - f(\eta,\xi)] \\
(4) \quad &+ \sum_{\substack{\eta(x)=\xi(x)=1 \\ \eta(z)=\xi(y)=0}} q(x,y,z,\xi,\eta)[f(\eta^{xz},\xi^{xy}) - f(\eta,\xi)] \\
(5) \quad &+ \sum_{\substack{\eta(x)=\xi(x)=1 \\ \eta(y)=0}} q(x,y,\eta,\xi)[f(\eta^{xy},\xi) - f(\eta,\xi)] \\
(6) \quad &+ \sum_{\substack{\eta(x)=\xi(x)=1 \\ \xi(y)=0}} q(x,y,\xi,\eta)[f(\eta,\xi^{xy}) - f(\eta,\xi)] \\
(7) \quad &+ \sum_{\substack{\eta(x)=1 \\ \eta(y)=\xi(x)=0}} q(x,y,\eta)[f(\eta^{xy},\xi) - f(\eta,\xi)] \\
(8) \quad &+ \sum_{\substack{\xi(x)=1 \\ \eta(x)=\xi(y)=0}} q(x,y,\xi)[f(\eta,\xi^{xy}) - f(\eta,\xi)] \\
(9) \quad &+ \sum_{\substack{\eta(x)=\xi(x)=\xi(y)=1 \\ \eta(y)=0}} q\delta(x,y,\eta,\xi)[f(\eta^{xy},\xi_x) - f(\eta,\xi)] \\
(10) \quad &+ \sum_{\substack{\eta(x)=\xi(x)=\eta(y)=1 \\ \xi(y)=0}} q\delta(x,y,\xi,\eta)[f(\eta_x,\xi^{xy}) - f(\eta,\xi)] \\
(11) \quad &+ \sum_{\eta(x)=1,\xi(x)=0} \delta(x,\eta)[f(\eta_x,\xi) - f(\eta,\xi)] \\
(12) \quad &+ \sum_{\eta(x)=0,\xi(x)=1} \delta(x,\xi)[f(\eta,\xi_x) - f(\eta,\xi)] \\
(13) \quad &+ \sum_{\eta(x)=1,\xi(x)=1} \delta(x,\eta\xi)[f(\eta_x,\xi_x) - f(\eta,\xi)],
\end{align}
$$

where

$$q(x,y,\eta\xi) = \mathbb{E}^x\left[\prod_{n=1}^{\sigma(y)-1} \eta(X_n)\xi(X_n), \sigma(y) < \sigma(x), \sigma(y) < \infty\right]$$

is the rate at which $\eta$ and $\xi$ particles at $x$ move together to $y \neq x$;

$q(x,y,z,\eta,\xi)$
$$= \mathbb{E}^x\left[\prod_{n=1}^{\sigma(y)-1} \eta(X_n) \times \prod_{n=1}^{\sigma(z)-1} \xi(X_n), \sigma(y) < \sigma(z), \sigma(z) < \sigma(x), \sigma(z) < \infty\right]$$

is the rate at which $\eta$ and $\xi$ particles at $x$ move, the $\eta$ particle stopping at $y \neq x$ and the $\xi$ particle continuing to $z \neq x$; $q(x,y,z,\xi,\eta)$ is as $q(x,y,z,\eta,\xi)$ with the roles of $\xi$ and $\eta$ interchanged;

$$q(x,y,\eta,\xi) = \mathbb{E}^x\left[\prod_{n=1}^{\sigma(y)-1} \eta(X_n) \times \prod_{n=1}^{\sigma(x)-1} \xi(X_n), \sigma(y) < \sigma(x) < \infty\right]$$

is the rate at which $\eta$ and $\xi$ particles at $x$ move to $y \neq x$ and remain at $x$, respectively; $q(x,y,\xi,\eta)$ is as $q(x,y,\eta,\xi)$ with the roles of $\xi$ and $\eta$ interchanged;

$$q\delta(x,y,\eta,\xi) = \mathbb{E}^x\left[\prod_{n=1}^{\sigma(y)-1} \eta(X_n) \times \prod_{n=1}^{\infty} \xi(X_n), \sigma(y) < \infty, \sigma(x) = \infty\right]$$

is the rate at which $\eta$ and $\xi$ particles at $x$ move to $y \neq x$ and disappear, respectively; $q\delta(x,y,\xi,\eta)$ is as $q\delta(x,y,\eta,\xi)$ with the roles of $\xi$ and $\eta$ interchanged and

$$\delta(x,\eta\xi) = \mathbb{E}^x\left[\prod_{n=1}^{\infty} \eta(X_n)\xi(X_n), \sigma(x) = \infty\right]$$

is the rate at which $\eta$ and $\xi$ particles at $x$ disappear; finally $q(x,y,\eta)$, $q(x,y,\xi)$, $\delta(x,\eta)$ and $\delta(x,\xi)$ are as in (1). Observe that

$$q(x,y,z,\eta,\xi) = q(x,y,\eta,\xi) = q\delta(x,y,\eta,\xi) = 0$$

unless $\xi(y) = 1$.

THEOREM 1.2. *Under the same hypothesis as in Theorem* 1.1, *let $\widetilde{\mu}$ be an invariant measure for the coupled long range exclusion process whose marginals $\mu_1$ and $\mu_2$ put no mass on $\mathbf{1}$. Then, for any cylinder function $f$, $\widetilde{L}f(\eta,\xi)$ is well defined $\widetilde{\mu}$ almost everywhere, belongs to $\mathbf{L}^1(\widetilde{\mu})$ and*

$$\int \widetilde{L}f \, d\widetilde{\mu} = 0.$$

REMARK. The proof of Theorem 1.2 will be omitted since it is similar to the one of Theorem 1.1. One just has to remark that $f$ is a linear combination of functions of the form $\prod_{(x,y)\in R_1 \times R_2} \eta(x)\xi(y)$, where $R_1$ and $R_2$ are finite subsets of $\mathbf{S}$ and prove the result for these functions proceeding as in the proof of our previous theorem.

Let $(\nu_\rho)_{\rho \in [0,1]}$ be the one-parameter family of Bernoulli product measures such that $\nu_\rho(\eta(x) = 1) = \rho$ for any $x \in \mathbf{S} = \mathbb{Z}^d$. Liggett has shown that, when $p(\cdot, \cdot)$ is a *random walk on* $\mathbb{Z}^d$, these measures are invariant:

THEOREM 1.3 (Theorem 4.2 and Corollary 4.4 of [6]). *Suppose* $p(y - x) = p(x, y)$ *for all* $x, y \in \mathbb{Z}^d$; *then* $\nu_\rho \in \mathcal{I}$ *for any* $\rho \in [0,1]$ *and in particular* $\nu_1 \in \mathcal{I}$.

In the last two sections of this paper, we apply Theorem 1.2 and coupling methods similar to [5] to show that in some cases all the invariant measures are convex combinations of the $\nu_\rho$'s ($0 \leq \rho \leq 1$).

The following two theorems require that $\mathbf{S}$ is an integer lattice and $p(x, y)$ is the transition matrix of an irreducible random walk on that lattice. We recall that by de Finetti's theorem (see Chapter VII.4 of [2]) the set of exchangeable probability measures coincides with the closed convex hull of $\{\nu_\rho : 0 \leq \rho \leq 1\}$.

THEOREM 1.4. *Suppose* $p(x, y)$ *is the transition matrix of an irreducible random walk on* $\mathbb{Z}^d$; *then the set of invariant and translation-invariant measures coincides with the set of exchangeable measures.*

This result is not new (see [3]). However, thanks to Theorem 1.2 we can give a more natural proof which shortcuts Guiol's.

THEOREM 1.5. *Suppose* $\mathbf{X} = \{0,1\}^\mathbb{Z}$, $p(x, y)$ *is the transition matrix of an irreducible random walk on* $\mathbb{Z}$ *such that* $\sum_x |x| p(x) < \infty$ *and* $\sum_x x p(x) = 0$. *Then*

$$\mathcal{I} = \{all\ exchangeable\ measures\}.$$

THEOREM 1.6. *Suppose* $\mathbf{X} = \{0,1\}^\mathbb{Z}$, $p(x, x+1) := p = 1 - p(x, x-1)$ *for all* $x \in \mathbb{Z}$. *Then*

$$\mathcal{I} = \{all\ exchangeable\ measures\}.$$

REMARK. Unlike the simple exclusion process, all the invariant measures of the one-dimensional nearest-neighbor long range exclusion process are translation invariant.

The paper is composed in the following way. Section 2 treats the construction of the process. In Section 3 we prove Theorem 1.1 and in Section 4 we give results relating the invariant measures of the coupled process to the invariant measures of the single process. The one-dimensional zero-mean case is treated in Section 5 where we also sketch the proof of Theorem 1.4, and Section 6 is devoted to the nearest-neighbor case. Throughout this paper $\mathbb{N}$ denotes the set of natural integers, that is, $\mathbb{N} = \{1, 2, \dots\}$; $\mathbb{Z}^+$ denotes the set of nonnegative integers, that is, $\mathbb{Z}^+ = \{0, 1, 2, \dots\}$ and $\mathrm{Card}(A)$ denotes the cardinality of $A$. Finally for $x, y \in \mathbb{R}$ we denote by $x \wedge y$ the minimum between $x$ and $y$.

**2. The construction.** In this section we construct the long range exclusion process simultaneously for all possible initial configurations. This construction follows [6] with a slightly different approach.

Let $(N_x(t))_{t \geq 0}$ $(x \in \mathbf{S})$ be a collection of Poisson processes of parameter 1 and for each $(n, x) \in \mathbb{N} \times \mathbf{S}$ let $(X_k^{n,x})_{k \in \mathbb{Z}^+}$ be a Markov chain with transition matrix $p(x, y)$ such that $X_0^{n,x} \equiv x$. Assume that all these Poisson processes and Markov chains are defined on the same probability space $\Omega$ and are independent. Let $\mathbf{X}_f$ be the set of finite initial configurations:

$$\mathbf{X}_f = \left\{ \eta \in \mathbf{X} : \sum_{x \in \mathbf{S}} \eta(x) < \infty \right\}.$$

We now provide a pathwise construction of the process starting from any $\eta \in \mathbf{X}_f$: A particle at $x \in \mathbf{S}$ will wait until the Poisson process $N_x(t)$ jumps. If this jump is the $n$th jump of $N_x(t)$, then the particle jumps to $X_{\tau(\xi)}^{n,x}$ where $\xi$ is the configuration of the process just before the jump and

$$\tau(\xi) = \inf\{r \geq 1 : X_r^{n,x} = x \text{ or } \xi(X_r^{n,x}) = 0\}.$$

Call $\eta_s^\xi$ the random configuration obtained at time $s$ when the initial configuration was $\xi$.

It is now easy to verify that there exists a subset $\Omega_0$ of $\Omega$ of probability 1 on which for any initial $\xi \in \mathbf{X}_f$:

(1) there are no simultaneous jumps of particles,
(2) the total number of jumps of $\eta_s^\xi$ is finite on any finite time interval,
(3) $\eta_0^\xi \equiv \xi$.

To extend the construction to infinite configurations we will need the following partial order in $\mathbf{X}$: we say that $\eta \leq \xi$ if $\eta(x) \leq \xi(x)$, $\forall x \in \mathbf{S}$. It is also easy to verify that our construction satisfies the following property:

(4) if $\zeta, \xi \in \mathbf{X}_f$ are such that $\zeta \leq \xi$, then on $\Omega_0$ we have $\eta_s^\zeta \leq \eta_s^\xi, \forall s \geq 0$.

For an arbitrary initial configuration $\xi \in \mathbf{X}$ and $s \geq 0$ we define on $\Omega_0$:

$$\eta_s^\xi(x) = \lim_{\zeta \uparrow \xi, \zeta \in \mathbf{X}_f} \eta_s^\zeta(x).$$

With this construction property (4) remains true for $\zeta, \xi \in \mathbf{X}$.

We have now defined on the same probability space $\eta_s^\xi$ for all $s \geq 0$ and all $\xi \in \mathbf{X}$. The coupled process starting from $(\zeta, \xi)$ is now defined as the process $(\eta_s^\zeta, \eta_s^\xi)$. When there will be no ambiguity on the initial configuration we will just denote by $\eta_t$ the process.

In the rest of this paper we denote by $\mathbf{P}$ and $\mathbf{E}$ the probability and the expectation operator of the space in which the process has been constructed, that is, the space $\Omega_0$. When computations involve auxiliary Markov chains or random variables, we call $\mathbb{P}$ and $\mathbb{E}$ the probability and the expectation operator which apply to them. Finally, when we need to emphasize the initial condition of one of these chains we will write $\mathbb{P}^x$ and $\mathbb{E}^x$.

**3. Invariant measures and generator.** This section is devoted to the proof of Theorem 1.1 which is performed in several lemmas. We begin by an informal description of the adopted strategy: Suppose $\xi_0$ is a finite initial configuration. Then standard results for Markov processes on countable state spaces yield

$$\mathbf{E}f(\xi_t) - f(\xi_0) = \mathbf{E}\int_0^t Lf(\xi_s)\,ds.$$

Then, if we are able to pass to the limit as $\xi_0$ grows to an infinite configuration and to integrate with respect to an invariant measure $\mu$ we should obtain

$$\int \left(\mathbf{E}\int_0^t Lf(\xi_s)\,ds\right)d\mu = 0.$$

Then, by Fubini and the invariance of $\mu$ we should have $t\int Lf(\xi_0)\,d\mu = 0$ as required. To justify this, certain integrability results are needed. These are provided by Lemmas 3.2 and 3.6.

3.1. *Notation and definitions.* Let $R$ be a finite subset of $\mathbf{S}$ and let $f_R(\eta) = \prod_{x\in R}\eta(x)$. When $R = \{x\}$ we will write $f_x$ for $f_R$. We decompose $L$ in two pieces $L = L^+ - L^-$ where

$$L^+ f_R(\eta) = \sum_{x\in R^c}\sum_{y\in R}\eta(x)q(x,y,\eta)[1-\eta(y)]\prod_{z\in R\setminus\{y\}}\eta(z)$$

and

$$L^- f_R(\eta) = \sum_{x\in R} f_R(\eta)\left(\delta(x,\eta) + \sum_{y\neq x}q(x,y,\eta)[1-\eta(y)]\right).$$

Note that $L^- f_R$ is bounded by the cardinal of $R$. Therefore the main difficulties to implement the strategy described at the beginning of this section are due to $L^+ f_R$. Finally let

$$L^{++} f_x(\eta) = \sum_{y\in\mathbf{S}}\eta(y)q(y,x,\eta).$$

This represents the rate at which particles jump to $x$ regardless of whether they may or may not stay at that site.

3.2. *Preliminary results.* The first lemma is a general result about Markov jump processes.

LEMMA 3.1. *Let $Q$ be the intensity matrix of a nonexplosive Markov jump process $X_t$ on a countable space $E$. Let $A$ be a nonempty subset of $E$, suppose $x \in E \setminus A$ and define $Q(x, A) = \sum_{y \in A} Q(x, y)$. Then for all $t \in \mathbb{R}^+$ we have*

$$\mathbb{P}^x\left(\int_0^{\tau_A} Q(X_s, A)\,ds \geq t\right) \leq e^{-t}$$

*where $\tau_A = \inf\{s > 0 : X_s \in A\}$.*

PROOF. Let

$$Q(x, A^c) = \sum_{y \in A^c \setminus \{x\}} Q(x, y) \quad \text{and} \quad Q(x) = Q(x, A) + Q(x, A^c).$$

Thus, $Q(x)$ is the rate at which the process leaves $x$. Let $\tau_0 = 0$, and for all $n \geq 1$ let $\tau_n = \inf\{t > \tau_{n-1} : X_t \neq X_{\tau_{n-1}}\}$ be the time of occurrence of the $n$th jump of the process. Then, let

$$F_x^n(t) = \mathbb{P}^x\left(\int_0^{\tau_A \wedge \tau_n} Q(X_s, A)\,ds > t\right)e^t.$$

To prove the statement we show by induction that $F_x^n(t) \leq 1$ for all $n \geq 1$, all $x \in A^c$ and all $t > 0$.

As $\tau_1 \leq \tau_A$ a direct computation gives

$$F_x^1(t) = \mathbb{P}^x\left(\int_0^{\tau_1} Q(X_s, A)\,ds > t\right)e^t = e^{-(Q(x)/Q(x,A))t}e^t = e^{-(Q(x,A^c)/Q(x,A))t}.$$

The second equality is valid since $X_s = x$ for all $0 \leq s < \tau_1$. Hence, $F_x^1(t) \leq 1$ for all $x \in E$, $A \subset E \setminus \{x\}$ and $t > 0$.

Now supposing we have the property up to rank $n$, we write

$$F_x^{n+1}(t) = \mathbb{P}^x(Q(x, A)\tau_1 > t)e^t$$
$$+ \mathbb{P}^x\left(Q(x, A)\tau_1 \leq t, \int_{\tau_1}^{\tau_A \wedge \tau_{n+1}} Q(X_s, A)\,ds > t - \tau_1 Q(x, A)\right)e^t.$$

The first term on the right-hand side is $F_x^1(t)$ and therefore equal to $e^{-(Q(x,A^c)/Q(x,A))t}$. Applying the strong Markov property to the second term we obtain

$$F_x^{n+1}(t) = e^{-(Q(x,A^c)/Q(x,A))t}$$
$$+ e^t\left[\int_0^{t/Q(x,A)} Q(x)e^{-Q(x)s}\right.$$

$$\times \sum_{y \in A^c \setminus \{x\}} \frac{Q(x,y)}{Q(x)} \mathbb{P}^y \left( \int_0^{\tau_A \wedge \tau_n} Q(X_u, A) \, du > t - sQ(x,A) \right) ds \Bigg]$$

$$= e^{-(Q(x,A^c)/Q(x,A))t}$$
$$+ \int_0^{t/Q(x,A)} e^{-Q(x)s} \sum_{y \in A^c \setminus \{x\}} Q(x,y) F_y^n(t - sQ(x,A)) e^{Q(x,A)s} \, ds$$

$$= e^{-(Q(x,A^c)/Q(x,A))t}$$
$$+ \sum_{y \in A^c \setminus \{x\}} Q(x,y) \int_0^{t/Q(x,A)} e^{-Q(x,A^c)s} F_y^n(t - sQ(x,A)) \, ds.$$

By the inductive hypothesis this implies

$$F_x^{n+1}(t) \leq e^{-(Q(x,A^c)/Q(x,A))t} + \int_0^{t/Q(x,A)} Q(x, A^c) e^{-Q(x,A^c)s} \, ds = 1$$

which concludes the proof of Lemma 3.1. □

The second result shows that an invariant measure, which does not charge the full configuration, gives weight 1 to the set of configurations for which the mean over time of the arrival rate of particles on a given site is finite.

LEMMA 3.2. *Suppose $p(\cdot, \cdot)$ is such that $\nu_1 \in \mathcal{I}$. Let $\mu \in \mathcal{I}$ such that $\mu(\{\mathbf{1}\}) = 0$ and let*

$$\mathcal{Y} = \left\{ \xi \in \mathbf{X} : \int_0^1 L^{++} f_x(\eta_s^\xi) \, ds < \infty \; \forall \, x \in \mathbf{S} \text{ a.s.} \right\}.$$

*Then $\mu(\mathcal{Y}) = 1$.*

The main idea in the proof of this lemma is the following: If $\int_0^1 L^{++} f_x(\eta_s^\xi) \, ds = \infty$, then infinitely many particles will attempt to move to $x$ in an arbitrarily small time interval. These particles will fill arbitrarily large neighborhoods of $x$ and the full configuration will be attained, contradicting the hypothesis $\mu(\{\mathbf{1}\}) = 0$. To develop this idea, we will need to introduce some notation and to prove an intermediate result (Lemma 3.3) which itself will be a consequence of a subsequent lemma.

Fix $x$ and let $\mathbf{S}_k$ be an increasing sequence of finite subsets of $\mathbf{S}$ such that:

(i) $x \in \mathbf{S}_1$,
(ii) for each $k$ and $y \in \mathbf{S}_k$ the Markov chain with transition matrix $p(x,y)$ starting at $x$ has a positive probability of hitting $y$ before leaving $\mathbf{S}_k$, and
(iii) $\bigcup_{k=1}^\infty \mathbf{S}_k = \mathbf{S}$.

For any arbitrary initial configuration $\xi \in \mathbf{X}$, and all $k \geq 1$, let $\xi^k(y) = \xi(y)\mathbf{I}_{\mathbf{S}_k}(y)$, where $\mathbf{I}_A(x) = 1$ if $x \in A$, 0 otherwise. Observe that for each finite $n$ the long range exclusion restricted to $\{\eta \in \mathbf{X} : \sum_x \eta(x) = n\}$ is a nonexplosive Markov chain on a countable state space. Let

$$\delta(\xi) = \mathbf{P}\left(\int_0^1 L^{++} f_x(\eta_s^\xi)\,ds = \infty\right)$$

and suppose $\ell$ is a strictly positive integer. Let $\tau_0 = 0$ and let $\tau_i = \tau_i(\ell)$, $i \geq 1$, be the successive times at which the Poisson process $\sum_{y \in \mathbf{S}_\ell}(N_y(t))_{t \geq 0}$, jumps.

Let

$$A(\xi) = \left\{\int_0^1 L^{++} f_x(\eta_s^\xi)\,ds = \infty\right\}$$

and for $j \geq 0$ let

$$A_j(\xi) = \left\{\int_{\tau_j \wedge 1}^{\tau_{j+1} \wedge 1} L^{++} f_x(\eta_s^\xi)\,ds = \infty\right\}.$$

For $k, \ell, j \in \mathbb{Z}^+$, let

$$C(\ell) = \{\eta \in \mathbf{X} : \eta(z) = 1 \ \forall z \in \mathbf{S}_\ell\},$$

$$\sigma^k(\ell, \xi) = \inf\{t \geq 0 : \eta_t^{\xi^k} \in C(\ell)\}$$

and

$$\sigma^k(\ell, j, \xi) = \inf\{t \geq \tau_j : \eta_t^{\xi^k} \in C(\ell)\}.$$

We will need the following lemma.

LEMMA 3.3. *With the preceding notation*

$$\lim_{k \to \infty} \mathbf{I}_{\{\sigma^k(\ell, j, \xi) < \tau_{j+1} \wedge 1\}} \geq \mathbf{I}_{A_j(\xi)} \qquad a.s.$$

PROOF OF LEMMA 3.2. For $k > \ell > n > 0$

$$\mathbf{P}(\eta_1^{\xi^k} \in C(n)) \geq \mathbf{P}(\eta_1^{\xi^k} \in C(n) | \sigma^k(\ell, \xi) < 1)\mathbf{P}(\sigma^k(\ell, \xi) < 1)$$

$$\geq \inf_{0 \leq s \leq 1} \mathbf{P}(\eta_s^{\mathbf{1}^\ell} \in C(n))\mathbf{P}(\sigma^k(\ell, \xi) < 1),$$

where the last inequality follows from the strong Markov property applied at time $\sigma^k(\ell, \xi)$ and the fact that for any $s \in [0, 1]$, $\mathbf{P}(\eta_s^\xi \in C(n)) \geq \mathbf{P}(\eta_s^\zeta \in C(n))$ if $\xi \geq \zeta$.

As $\sigma^k(\ell, \xi) \leq \sigma^k(\ell, j, \xi)$ for all $j \geq 0$, Lemma 3.3 implies that

$$\lim_{k \to \infty} \mathbf{I}_{\{\sigma^k(\ell, \xi) < 1\}} \geq \mathbf{I}_{\bigcup_j A_j(\xi)} = \mathbf{I}_{A(\xi)} \qquad a.s.$$

Therefore
$$\lim_{k\to\infty} \mathbf{P}(\eta_1^{\xi^k} \in C(n)) \geq \inf_{0\leq s\leq 1} \mathbf{P}(\eta_s^{\mathbf{1}^\ell} \in C(n))\mathbf{P}(A(\xi)).$$

Hence
$$\mathbf{P}(\eta_1^{\xi} \in C(n)) \geq \mathbf{P}(A(\xi))\lim_{\ell\to\infty}\inf_{0\leq s\leq 1}\mathbf{P}(\eta_s^{\mathbf{1}^\ell} \in C(n)) = \mathbf{P}(A(\xi)) = \delta(\xi),$$

where the first equality in the last formula is due to the invariance of the full configuration. Hence $\mathbf{P}(\eta_1^\xi \equiv 1) = \lim_n \mathbf{P}(\eta_1^\xi \in C(n)) \geq \delta(\xi)$, and
$$0 = \mu(\{\mathbf{1}\}) = \mu S(1)(\{\mathbf{1}\}) \geq \int \delta(\xi)\, d\mu(\xi).$$

Hence $\delta(\xi) = 0$ $\mu$-a.e. This concludes the proof of Lemma 3.2. □

It remains to prove Lemma 3.3. To do so, we introduce some extra notation: let $E_0 = \{x\}$, for $i \geq 1$ we denote by $E_i$ the set of points in $\mathbf{S}_\ell$ reachable by $X_n$ visiting $i-1$ intermediate points in $\mathbf{S}_\ell$, that is,
$$E_i = E_i(\ell) = \left\{ y \in \mathbf{S}_\ell \colon \sum_{x_1,\ldots,x_{i-1}\in\mathbf{S}_\ell} p(x,x_1)\cdots p(x_{i-1},y) > 0 \right\}.$$

For $\zeta \in \mathbf{X}_f$ let
$$F_0 = F_0(\zeta,\ell) = E_0 \cap \{y \in \mathbf{S}_\ell : \zeta(y) = 0\},$$

and for all $i \geq 1$ let
$$F_i = F_i(\zeta,\ell) = E_i \cap \{y \in \mathbf{S}_\ell : \zeta(y) = 0\} \setminus \bigcup_{j=0}^{i-1} F_j.$$

Since $\mathbf{S}_\ell$ is finite, there exists $r < \infty$ such that
$$\bigcup_{j=0}^{r} F_j = \{y \in \mathbf{S}_\ell : \zeta(y) = 0\}.$$

Let
$$b(\ell) = \inf_{i\leq r}\inf_{y\in E_i} \sum_{x_j \in E_j, j=1,\ldots,i-1} p(x,x_1)\cdots p(x_{i-1},y),$$

and note that $0 < b(\ell) \leq 1$. Now we state and prove another lemma and then we derive Lemma 3.3 from it.

LEMMA 3.4. *Fix $\varepsilon > 0$ and $c \in (0,1]$, let $m_0 = \operatorname{Card}(\mathbf{S}_\ell)$ and let $a$ be a positive real number such that $e^{-a} < \varepsilon/m_0$. Suppose $\zeta \in \mathbf{X}_f$ and let*
$$B_c(\zeta) = \left\{ \int_0^{\tau_1 \wedge c} b(\ell) L^{++} f_x(\eta_s^\zeta)\, ds \geq am_0 \right\},$$

*then*

$$\mathbf{P}(B_c(\zeta) \cap \{\sigma(\ell,\zeta) \geq \tau_1 \wedge c\}) \leq \varepsilon,$$

*where*

$$\sigma(\ell,\zeta) = \inf\{t \geq 0 : \eta_t^\zeta \in C(\ell)\}.$$

PROOF. Enumerate the points of $\bigcup_{i=0}^r F_i$ starting with points in $F_0$ (if any), then points in $F_1$, and so on. Let $y_1, \ldots, y_m$ be this enumeration. Define

$$D_{y_j} = \{\eta : \eta(y_j) = 1\},$$
$$\rho_0 = 0,$$
$$\rho_1 = \inf\{t > 0 : \eta_t^\zeta(y_1) = 1\},$$
$$\rho_i = \inf\{t \geq \rho_{i-1} : \eta_t^\zeta(y_i) = 1\}, \qquad i = 2, \ldots, m.$$

Since $m \leq m_0$, we have

$$\mathbf{P}\bigg(\int_0^{\rho_m \wedge \tau_1} b(\ell) L^{++} f_x(\eta_s^\zeta) \, ds \geq a m_0\bigg)$$

(14) $$\leq \mathbf{P}\bigg(\int_0^{\rho_1 \wedge \tau_1} b(\ell) L^{++} f_x(\eta_s^\zeta) \, ds \geq a\bigg)$$

(15) $$+ \sum_{i=1}^{m-1} \mathbf{P}\bigg(\int_{\rho_i \wedge \tau_1}^{\rho_{i+1} \wedge \tau_1} b(\ell) L^{++} f_x(\eta_s^\zeta) \, ds \geq a\bigg).$$

Since the process $\eta_s^\zeta$ is a nonexplosive Markov jump process on a countable state space, we may apply to it Lemma 3.1. Using the notation of that lemma, $s \in (0, \rho_1 \wedge \tau_1)$ implies that

$$b(\ell) L^{++} f_x(\eta_s^\zeta) \leq Q(\eta_s^\zeta, D_{y_1}),$$

where $Q$ is the intensity matrix of $\eta_s^\zeta$.

Therefore, it follows from Lemma 3.1 that (14) is bounded above by

$$\mathbf{P}\bigg(\int_0^{\rho_1} Q(\eta_s^\zeta, D_{y_1}) \, ds \geq a\bigg) \leq \frac{\varepsilon}{m_0}.$$

Similarly, for $1 \leq i \leq m-1$ each corresponding term in (15) is bounded above by $\varepsilon/m_0$. Therefore

$$\mathbf{P}\bigg(\int_0^{\rho_m \wedge \tau_1} b(\ell) L^{++} f_x(\eta_s^\zeta) \, ds \geq a m_0\bigg) \leq \varepsilon \frac{m}{m_0} \leq \varepsilon.$$

Hence

$$\mathbf{P}\bigg(\int_0^{\rho_m \wedge \tau_1} b(\ell) L^{++} f_x(\eta_s^\zeta) \, ds < a m_0, B_c(\zeta)\bigg) \geq \mathbf{P}(B_c(\zeta)) - \varepsilon.$$

Since on $B_c(\zeta)$

$$am_0 \leq \int_0^{\tau_1 \wedge c} b(\ell) L^{++} f_x(\eta_s^\zeta)\, ds,$$

this implies that

$$\mathbf{P}(\rho_m < \tau_1 \wedge c, B_c(\zeta)) \geq \mathbf{P}(B_c(\zeta)) - \varepsilon.$$

It now follows from

$$\{\rho_m < \tau_1 \wedge c, B_c(\zeta)\} \subset \{\eta_{\rho_m}^\zeta \in C(\ell), B_c(\zeta)\},$$

that

$$\mathbf{P}(\sigma(\ell,\zeta) < \tau_1 \wedge c, B_c(\zeta)) \geq \mathbf{P}(B(\zeta)) - \varepsilon,$$

which implies the lemma. $\square$

PROOF OF LEMMA 3.3. Fix $\varepsilon > 0$, let $a$ be as in Lemma 3.4 and let

$$A_j(\xi, k) = \left\{ \int_{\tau_j \wedge 1}^{\tau_{j+1} \wedge 1} b(\ell) L^{++} f_x(\eta_s^{\xi^k})\, ds \geq am_0 \right\}.$$

Applying the strong Markov property to the process $\eta_s^{\xi^k}$ at the stopping time $\tau_j$ we get

$$\mathbf{P}(A_j(\xi, k) \cap \{\sigma^k(\ell, j, \xi) \geq \tau_{j+1} \wedge 1\})$$
$$= \int \mathbf{P}(B_c(\zeta) \cap \{\sigma(\ell, \zeta) \geq \tau_1 \wedge c\})\, d\mu(\zeta, c),$$

where $\mu(\zeta, c)$ is the joint distribution of $(\eta_{\tau_j}^{\xi^k}, 1 - \tau_j)$ restricted to $\{\tau_j < 1\}$. It now follows from Lemma 3.4 that the right-hand side of this equation is less than $\varepsilon$. Since by monotone convergence

$$\lim_{k \to \infty} \int_{\tau_j \wedge 1}^{\tau_{j+1} \wedge 1} b(\ell) L^{++} f_x(\eta_s^{\xi^k})\, ds = \int_{\tau_j \wedge 1}^{\tau_{j+1} \wedge 1} b(\ell) L^{++} f_x(\eta_s^\xi)\, ds,$$

there exists $k_0$ such that

$$\mathbf{P}(A_j(\xi) \cap A_j(\xi, k)^c) \leq \varepsilon \qquad \forall k \geq k_0.$$

Hence,

$$\mathbf{P}(A_j(\xi) \cap \{\sigma^k(\ell, j, \xi) \geq \tau_{j+1} \wedge 1\}) \leq 2\varepsilon$$

holds for all $k \geq k_0$. Since $\varepsilon$ is arbitrary and $\{\sigma^k(\ell, j, \xi) \geq \tau_{j+1} \wedge 1\}$ decreases with $k$, the lemma follows. $\square$

LEMMA 3.5. *Let $R$ be a finite subset of $\mathbf{S}$. Then on $\Omega_0$*

$$\text{(16)} \qquad \lim_{k \to \infty} L^+ f_R(\eta_s^{\xi^k}) = L^+ f_R(\eta_s^\xi) \qquad \forall \xi \in \mathbf{X}, s \in \mathbb{R}^+,$$

and

$$\text{(17)} \qquad \lim_{k \to \infty} L^- f_R(\eta_s^{\xi^k}) = L^- f_R(\eta_s^\xi) \qquad \forall \xi \in \mathbf{X}, s \in \mathbb{R}^+,$$

PROOF. To prove the first equality we consider three cases:

(a) If for $x, y \in R$, $y \neq x$, $\eta_s^\xi(x) = \eta_s^\xi(y) = 0$, then $\eta_s^{\xi^k}(x) = \eta_s^{\xi^k}(y) = 0$ for all $k$ and $L^+ f_R(\eta_s^{\xi^k}) = L^+ f_R(\eta_s^\xi) = 0$.

(b) If $x$ is the only point in $R$ such that $\eta_s^\xi(x) = 0$, then $\eta_s^{\xi^k}(x) = 0$ for all $k$. Hence, the result follows from monotone convergence.

(c) If $f_R(\eta_s^\xi) = 1$, then $f_R(\eta_s^{\xi^k}) = 1$ for $k$ large enough and $\lim_{k \to \infty} L^+ f_R(\eta_s^{\xi^k}) = L^+ f_R(\eta_s^\xi) = 0$.

For the second equality we consider two cases:

(a') If $f_R(\eta_s^\xi) = 0$, then $f_R(\eta_s^{\xi^k}) = 0$ for all $k$. This implies that $L^- f_R(\eta_s^{\xi^k}) = 0 = L^- f_R(\eta_s^\xi)$ for all $k$.

(b') If $f_R(\eta_s^\xi) = 1$, then $f_R(\eta_s^{\xi^k}) = 1$ for $k$ large enough, and the result follows from the bounded convergence theorem. □

Let $Z_n$, $n \geq 1$, be a sequence of i.i.d. random variables, whose common distribution is the exponential of parameter 1, and let $N$ be a Poisson random variable of parameter 1 and independent of the $Z$'s. Then define $U = \sum_{i=1}^{N+1} Z_i$ and $U_n = \sum_{i=1}^{(N+1) \wedge n} Z_i$.

LEMMA 3.6. *Let $U$ be as above. Then,*

$$\text{(18)} \qquad \mathbf{P}\left(\int_0^1 L^+ f_x(\eta_s^\xi) \, ds > a\right) \leq \mathbb{P}(U > a) \qquad \forall \xi \in \mathbf{X}, a \in \mathbb{R}^+, x \in \mathbf{S}.$$

PROOF. We start proving the lemma for $\xi \in \mathbf{X}_f$. Define the following stopping times:

$$\sigma^0 = 0,$$
$$\tau^i = \inf\{t > \sigma^{i-1} : \eta_t^\xi(x) = 1\}, \qquad i \geq 1,$$
$$\sigma^i = \inf\{t > \tau^i : \eta_t^\xi(x) = 0\}, \qquad i \geq 1.$$

Let

$$X_i = \int_{\sigma^{i-1}}^{\tau^i} L^+ f_x(\eta_s^\xi) \, ds$$

and
$$Y_i = \sigma^i - \tau^i.$$
Then applying Lemma 3.1 to the set $A = \{\xi \in \mathbf{X}_f : \xi(x) = 1\}$ we get

(19) $\quad \mathbf{P}\left(\int_0^{\tau^1} L^+ f_x(\eta_s^\xi)\, ds \geq t\right) \leq e^{-t} \qquad \forall \xi \in \mathbf{X}_f, t \in \mathbb{R}^+.$

Then note that each time $x$ is occupied, it remains so for at least an exponential time of parameter 1. Therefore, the number of subintervals of $[0,1]$ during which $x$ is not occupied is stochastically dominated by $N+1$. Since during the time intervals in which $x$ is occupied $L^+ f_x(\eta) = 0$,

$$\int_0^1 L^+ f_x(\eta_s^\xi)\, ds \leq \sum_{i \geq 1} \mathbf{I}_{[0,1)}(\sigma^{i-1}) X_i \leq \sum_{i \geq 1} X_i \mathbf{I}_{[0,1)}\left(\sum_{j=1}^{i-1} Y_j\right),$$

where by convention $\sum_{j=1}^{i-1} Y_j = 0$ if $i = 1$.

Let $F_n(x_1, y_1, \ldots, x_n)$ be a bounded measurable function, increasing in the variables $x$ and decreasing in the variables $y$. We will now prove that

(20) $\quad \mathbf{E}(F_n(X_1, Y_1, \ldots, X_n)) \leq \mathbb{E}(F_n(Z_1, \ldots, Z_{2n-1})),$

where the random variables $Z_i$ are as above.

This is done by induction on $n$. For $n = 1$ this is a consequence of (19). Conditioning with respect to the $\sigma$-algebras $\mathcal{F}_{\sigma^{n-1}}$ we get

$$\mathbf{E}(F_n(X_1, Y_1, \ldots, X_n)|\mathcal{F}_{\sigma^{n-1}}) = \int F_n(X_1, Y_1, \ldots, Y_{n-1}, x_n)\, d\mathbf{P}_{X_n|\mathcal{F}_{\sigma^{n-1}}}(x_n).$$

Since $F_n$ is increasing in $x_n$ and by the strong Markov property and Lemma 3.1 the conditional distribution of $X_n$ given $\mathcal{F}_{\sigma^{n-1}}$ is a.s. stochastically bounded above by the distribution of $Z_{2n-1}$, the right-hand side is a.s. bounded above by

$$\int F_n(X_1, Y_1, \ldots, Y_{n-1}, x_n)\, d\mathbb{P}_{Z_{2n-1}}(x_n).$$

Conditioning this last expression on $\mathcal{F}_{\tau^{n-1}}$ and arguing in a similar way we get

$\mathbf{E}(F_n(X_1, Y_1, \ldots, X_n)|\mathcal{F}_{\tau^{n-1}})$
$$\leq \iint F_n(X_1, Y_1, \ldots, X_{n-1}, y_{n-1}, x_n)\, d\mathbb{P}_{Z_{2n-2}}(y_{n-1})\, d\mathbb{P}_{Z_{2n-1}}(x_n)$$

and (20) follows by the inductive hypothesis.

Applying (20) to
$$F_n(x_1, y_1, \ldots, x_n) = \begin{cases} 1, & \text{if } \sum_{i \geq 1}^n x_i \mathbf{I}_{[0,1)}\left(\sum_{j=1}^{i-1} y_j\right) > a, \\ 0, & \text{otherwise,} \end{cases}$$

we get
$$\mathbf{P}\left(\sum_{i\geq 1}^{n} X_i \mathbf{I}_{[0,1)}\left(\sum_{j=1}^{i-1} Y_j\right) > a\right) \leq \mathbb{P}(U_n > a),$$
and taking limits as $n$ goes to infinity we obtain (18).

To prove (18) for an infinite initial configuration take an increasing limit of finite initial configurations and deduce it from the finite case using (16) and Fatou's lemma. $\square$

3.3. *Proof of Theorem* 1.1. As $L^+ f_R(\eta_s^{\xi^k}) \leq \sum_{x\in R} L^+ f_x(\eta_s^{\xi^k}) \leq \sum_{x\in R} L^{++} f_x(\eta_s^{\xi})$, by (16), the dominated convergence theorem and the definition of $\mathcal{Y}$ we have
$$\lim_{k\to\infty} \int_0^1 L^+ f_R(\eta_s^{\xi^k}) \, ds = \int_0^1 L^+ f_R(\eta_s^{\xi}) \, ds \qquad \text{a.s.}, \forall \xi \in \mathcal{Y}.$$
As $|L^- f_R| \leq \mathrm{Card}(R)$ a similar equality for $L^-$ follows from (17). Therefore,

$$\text{(21)} \qquad \lim_{k\to\infty} \int_0^1 L f_R(\eta_s^{\xi^k}) \, ds = \int_0^1 L f_R(\eta_s^{\xi}) \, ds \qquad \text{a.s.}, \forall \xi \in \mathcal{Y}.$$

Since $L^+ f_R(\eta_s^{\xi^k}) \leq L^+ \sum_{x\in R} f_x(\eta_s^{\xi^k})$, and by (18) the sequence (in $k$) of random variables
$$\left\{\int_0^1 L^+ f_x(\eta_s^{\xi^k}) \, ds\right\}_{k\in\mathbb{N}}$$
is uniformly integrable, it follows from $|L^- f_R| \leq \mathrm{Card}(R)$ and (21) that

$$\text{(22)} \qquad \mathbf{E}\left[\int_0^1 L f_R(\eta_s^{\xi^k}) \, ds\right] \to \mathbf{E}\left[\int_0^1 L f_R(\eta_s^{\xi}) \, ds\right] \qquad \forall \xi \in \mathcal{Y} \text{ when } k \uparrow \infty.$$

As the process $\eta_s^{\xi^k}$ is a Markov process on the countable space $\mathbf{X}_f$ with a bounded $Q$ matrix we have
$$\mathbf{E} f_R(\eta_1^{\xi^k}) - f_R(\xi^k) = \mathbf{E}\int_0^1 L f_R(\eta_s^{\xi^k}) \, ds.$$
Since $|\mathbf{E} f_R(\eta_1^{\xi^k}) - f_R(\xi^k)| \leq 1$, it follows from (22) and the dominated convergence theorem that
$$\lim_{k\to\infty} \int (\mathbf{E} f_R(\eta_1^{\xi^k}) - f_R(\xi^k)) \, d\mu(\xi) = \int \mathbf{E}\left(\int_0^1 L f_R(\eta_s^{\xi}) \, ds\right) d\mu(\xi).$$
By the dominated convergence theorem and the invariance of $\mu$ the left-hand side above is 0. Since $|L^- f_R| \leq \mathrm{Card}(R)$ we may apply Tonelli's theorem to the right-hand side. It then follows from the invariance of $\mu$ that this right-hand side is equal to
$$\int L f_R(\eta) \, d\mu(\eta),$$
thus the theorem is proved.

**4. Coupled process and invariant measures.** In this section we give some general results relating the invariant measures of the single process to the invariant measures of the coupled process. Let $\mathcal{I}$ and $\widetilde{\mathcal{I}}$ be the set of invariant measures for the single process and the set of invariant measures for the coupled process, respectively, and let $\mathcal{S}$ and $\widetilde{\mathcal{S}}$ be the set of translation-invariant measures on $\mathbf{X}$ and $\mathbf{X}^2$, respectively.

LEMMA 4.1. *Suppose $\mu_1, \mu_2 \in \mathcal{I}$. Then there exists a probability measure $\widetilde{\nu}$ on $\mathbf{X} \times \mathbf{X}$ which is invariant for the coupled process and has marginals $\mu_1$ and $\mu_2$. Moreover, if $\mathbf{S} = \mathbb{Z}^d$, and $p(x,y)$, $\mu_1$ and $\mu_2$ are translation invariant, then we can choose a translation-invariant measure $\widetilde{\nu}$ as above.*

For a proof of this lemma we refer the reader to Proposition 5.4 and Remark 5.5 in [3].

LEMMA 4.2. *Let $\widetilde{\nu}$ be an invariant measure for the coupled process such that*

$$\widetilde{\nu}(\{(\eta, \xi) : \eta(x) = \xi(y) = 1, \eta(y) = \xi(x) = 0\}) = 0,$$

*whenever $p(x,y) > 0$. If $p(x,y)$ is irreducible, then*

$$\widetilde{\nu}(\{(\eta, \xi) : \eta \leq \xi \text{ or } \xi \leq \eta\}) = 1.$$

Thanks to Theorem 1.2 the proof of this lemma follows a standard induction argument (see Lemma 2.5 in [5]) and will be omitted.

LEMMA 4.3. *Suppose that $\mathbf{S}$ is the $d$-dimensional integer lattice $\mathbb{Z}^d$ and that for any $\mu \in \mathcal{I}$ and any $\rho \in [0,1]$ there exists a probability measure $\widetilde{\nu}_\rho$ on $\mathbf{X}^2$ such that:*

(i) *$\widetilde{\nu}_\rho$ is invariant for the coupled process,*
(ii) *its first marginal is $\nu_\rho$, its second marginal is $\mu$ and*
(iii) *$\widetilde{\nu}_\rho\{(\eta, \xi) : \eta \leq \xi \text{ or } \xi \leq \eta\} = 1$.*

*Then all elements of $\mathcal{I}$ are exchangeable.*

PROOF. For $\eta \in \mathbf{X}$ define

$$U(\eta) = \limsup_n \frac{1}{(2n+1)^d} \sum_{(x_1,\ldots,x_d) \in \mathbb{Z}^d : |x_i| \leq n, i=1,\ldots,d} \eta(x)$$

and

$$L(\eta) = \liminf_n \frac{1}{(2n+1)^d} \sum_{(x_1,\ldots,x_d) \in \mathbb{Z}^d : |x_i| \leq n, i=1,\ldots,d} \eta(x).$$

Let $\mu$ be an element of $\mathfrak{I}$ and let $\delta > 0$. Then, as in page 426 in [1] it follows from the hypothesis of this lemma that

$$U(\eta) = L(\eta), \qquad \mu\text{-a.e.}$$

Let $0 = a_0 < a_1 < \cdots < a_{r-1} < a_r = 1$ be such that $a_i - a_{i-1} < \delta$ for $i = 1, \ldots, r$ and $\mu(\{\eta : U(\eta) = a_i\}) = 0$ for $i = 1, \ldots, r-1$. Then, for $i = 1, \ldots, r$ let

$$A_i = \{\eta \in \mathbf{X} : U(\eta) = L(\eta) \in [a_{i-1}, a_i]\},$$
$$B_i = \{\eta \in \mathbf{X} : U(\eta) = L(\eta) < a_i\},$$

and for $i = 0, \ldots, r-1$ let

$$C_i = \{\eta \in \mathbf{X} : U(\eta) = L(\eta) > a_i\}.$$

We will now show that under the assumption $\mu(A_i) > 0$ the conditional measure $\mu(\cdot|A_i)$ is in $\mathfrak{I}$. This is done in two steps: First let $\widetilde{\nu}_i$ be an invariant measure for the coupled process whose marginals are $\mu$ and $\nu_{a_{i-1}}$ and concentrates on $\{(\eta, \xi) : \eta \geq \xi \text{ or } \eta \leq \xi\}$. Since the set $\{(\eta, \xi) : \eta \geq \xi\}$ is stable for the coupled process, it follows that the conditional measure $\widetilde{\nu}_i(\cdot|\{(\eta, \xi) : \eta \geq \xi\})$ is also invariant for the coupled process. Hence, its first marginal $\mu(\cdot|C_{i-1}) \in \mathfrak{I}$. Applying a similar argument to this last measure we see that $\mu(\cdot|B_i \cap C_{i-1}) = \mu(\cdot|A_i) \in \mathfrak{I}$. Finally applying the hypothesis of the proposition to $\mu(\cdot|A_i)$ we see that $\nu_{a_{i-1}} \leq \mu(\cdot|A_i) \leq \nu_{a_i}$. Therefore, there exist nonnegative reals $\lambda_1, \ldots, \lambda_r$ whose sum is 1 and satisfy

$$\sum_{i=1}^{r} \lambda_i \nu_{a_{i-1}} \leq \mu \leq \sum_{i=1}^{r} \lambda_i \nu_{a_i}.$$

From this double inequality we conclude that for any $n$ and any pair of subsets $A$ and $B$ of $\mathbb{Z}^d$ of cardinality $n$,

$$\left|\mu\left(\left\{\eta : \prod_{x \in A} \eta(x) = 1\right\}\right) - \mu\left(\left\{\eta : \prod_{x \in B} \eta(x) = 1\right\}\right)\right|$$

is bounded above by an expression which tends to 0 with $\delta$. Since $\delta$ is arbitrary it follows that $\mu(\{\eta : \prod_{x \in A} \eta(x) = 1\}) = \mu(\{\eta : \prod_{x \in B} \eta(x) = 1\})$. Hence $\mu$ is exchangeable. $\square$

REMARK. In view of the last two lemmas, to prove the last three theorems of this paper it suffices to show for any $\mu \in \mathfrak{I}$ ($\mu \in \mathfrak{I} \cap \mathcal{S}$) and any $\rho \in [0, 1]$, there exists an invariant (and translation-invariant) measure for the coupled process $\widetilde{\nu}_\rho$ with marginals $\mu$ and $\nu_\rho$ such that

$$\widetilde{\nu}_\rho(\{(\eta, \xi) : \eta(x) = \xi(y) = 1, \eta(y) = \xi(x) = 0\}) = 0,$$

whenever $p(x,y) > 0$. Moreover, we may assume that $\mu(\{\mathbf{1}\}) = 0$, since otherwise we can decompose $\mu$ as $\mu(\{\mathbf{1}\})\nu_1 + (1 - \mu(\{\mathbf{1}\}))\mu'$ where $\mu'$ is an invariant measure such that $\mu'(\{\mathbf{1}\}) = 0$. Furthermore, we may also assume that $\rho \in (0,1)$, since the hypothesis of Lemma 4.3 follows immediately from Lemma 4.1 if $\rho = 0$ or $1$.

**5. Random walk with zero mean.** As remarked at the end of the previous section, to prove Theorem 1.5 (Theorem 1.4) of this paper it suffices to show for any $\mu \in \mathcal{I}$ ($\mu \in \mathcal{I} \cap \mathcal{S}$) such that $\mu(\{\mathbf{1}\}) = 0$ and any $\rho \in (0,1)$, there exists an invariant (and translation-invariant) measure for the coupled process $\widetilde{\nu}_\rho$ with marginals $\nu_\rho$ and $\mu$ such that

$$(23) \quad \widetilde{\nu}_\rho(\{(\eta,\xi) : \eta(x) = \xi(y) = 1, \eta(y) = \xi(x) = 0\}) = 0 \qquad \text{if } p(x,y) > 0.$$

Since by Lemma 4.1 there is a measure $\widetilde{\nu}_\rho \in \widetilde{\mathcal{I}}$ ($\widetilde{\mathcal{I}} \cap \widetilde{\mathcal{S}}$) with marginals $\nu_\rho$ and $\mu$, we only have to prove that it satisfies (23).

In this section (except the very end of it where we sketch the proof of Theorem 1.4) we restrict ourselves to dimension 1: $\mathbf{X} = \{0,1\}^{\mathbb{Z}}$ and we suppose that the underlying transition matrix $p$ corresponds to a zero-mean random walk, that is, for all $x, y \in \mathbb{Z}$, $p(y - x) = p(x, y)$,

$$(24) \qquad \sum_{x \in \mathbb{Z}} |x| p(x) < \infty$$

and

$$(25) \qquad \sum_{x \in \mathbb{Z}} x p(x) = 0.$$

5.1. *Notation and definitions.* Recall that for $x \in \mathbb{Z}$, $\sigma(y)$ denotes the hitting time of $\{y\}$ for the random walk with transition matrix $p(x,y)$. We now change slightly the notation. For any $\eta \in \mathbf{X}$ we denote

$$\sigma(\eta) = \inf\{k \geq 1 : \eta(X_k) = 0\}$$

where $(X_n)_{n \in \mathbb{N}}$ is a random walk with transition probability $p$. Let $p^*$ be the transition probability of the reverse walk defined by $p^*(x,y) = p(y,x)$ for all $x, y \in \mathbb{Z}$. Observe that the Markov chain $X_k^*$ with transition matrix $p^*$ is also a zero-mean random walk. Now, $q^*(x, y, \eta)$ and $\sigma^*(\eta)$ are defined as $q$ and $\sigma$ substituting $X_k^*$ for $X_k$.

Let

$$(26) \qquad \overline{q}(x,y,\eta) = \mathbb{E}^x\left[\prod_{k=1}^{\sigma(y)-1} \eta(X_k), \sigma(y) < \infty\right].$$

Let $Y_0 = X_0$ and $Y_n = X_n - X_{n-1}$ for $n \geq 1$ be the increments of the walk. The $(Y_i)$'s for $i \geq 1$ are i.i.d. random variables with distribution $p$. For a fixed $\eta \in \mathbf{X}$ observe that $\sigma(\eta)$ is a stopping time for the natural filtration of the $(Y_i)$'s, and for all $x \neq y \in \mathbb{Z}$

$$(27) \qquad [1 - \eta(y)]\overline{q}(x, y, \eta) = \mathbb{P}^x\left(\sum_{n=0}^{\sigma(\eta)} Y_n = y, \sigma(\eta) < \infty\right).$$

The couple $(\eta_t, \xi_t)$ will designate the basic coupled process. We will say that there is a *positive [negative] discrepancy* at site $x$ if $\eta(x) > \xi(x)$ [$\eta(x) < \xi(x)$]. We now define:

$$A(x, y, \eta, \xi)$$

$$(28) \quad = \mathbf{I}_{\eta(x)=1,\xi(x)=0}\bigg(\mathbf{I}_{\eta(y)=\xi(y)=0}\, q(x, y, \eta)$$

$$+ \sum_{z \neq x} \mathbf{I}_{\eta(z)=\xi(z)=1,\eta_z(y)=\xi_z(y)=0}\, q(z, x, \eta\xi)\overline{q}(x, y, \eta_z)\bigg).$$

This represents the rate at which a positive discrepancy at site $x$ reaches a vacant site $y$: The term preceding the "+" sign considers the case in which the movement of the $\eta$ particle starts at $x$ while the term after that sign takes into account a movement due to the simultaneous arrival from $z$ to $x$ of $\eta$ and $\xi$ particles. We define also:

$$B(x, y, \eta, \xi)$$

$$(29) \quad = \mathbf{I}_{\eta(x)=1,\xi(x)=0}\bigg(\mathbf{I}_{\eta(y)=0}\, q(x, y, \eta)$$

$$+ \sum_{z \neq x} \mathbf{I}_{\eta(z)=\xi(z)=1,\eta_z(y)=0}\, q(z, x, \eta\xi)\overline{q}(x, y, \eta_z)\bigg).$$

This is almost the same expression as (28) but site $y$ may or may not be occupied on the $\xi$ coordinate; thus the positive discrepancy may or may not coalesce with a negative one. We also define:

$$(30) \quad C(x, y, \eta) = \mathbf{I}_{\eta(x)=1}\bigg(\mathbf{1}_{\eta(y)=0}\, q(x, y, \eta)$$

$$+ \sum_{z \neq x} \mathbf{I}_{\eta(z)=1,\eta_z(y)=0}\, q(z, x, \eta)\overline{q}(x, y, \eta_z)\bigg).$$

This is like the previous expression (29) without any restriction on the $\xi$ coordinate.

Finally let

$$(31) \qquad D(x, y, \eta, \xi) = q(x, y, \eta)\mathbf{1}_{\eta(x)=\xi(y)=1,\xi(x)=\eta(y)=0},$$

which is the rate at which a positive discrepancy at site $x$ moves and coalesces with a negative one at site $y$.

REMARK. Obviously for all $\eta, \xi \in \mathbf{X}$ and $x, y \in \mathbb{Z}$
$$A(x,y,\eta,\xi) \leq B(x,y,\eta,\xi) \leq C(x,y,\eta).$$

As before $\nu_\rho$ denotes the Bernoulli product measure with constant density $\rho \in (0,1)$. Let $\mu \in \mathcal{I}$ and denote by $\widetilde{\nu}_\rho$ an invariant measure for the coupled process with marginals $\nu_\rho$ and $\mu$.

Let $A_\rho(x,y) = \int A(x,y,\eta,\xi) \, d\widetilde{\nu}_\rho$. Similarly, let $B_\rho(x,y)$, $C_\rho(x,y)$ and $D_\rho(x,y)$ be the respective integrals of (29), (30) and (31) with respect to $\widetilde{\nu}_\rho$. From our previous remark we get

$$(32) \qquad A_\rho(x,y) \leq B_\rho(x,y) \leq C_\rho(x,y).$$

The proof of Theorem 1.5 relies on the following observations:

(1) jumps of positive discrepancies have mean 0,
(2) when discrepancies of different sign meet they disappear,
(3) the number of positive discrepancies in the interval $[-n,n]$ can only increase through boundary effects, but due to (1) the Cesàro averages on $n$ of these effects tend to 0,
(4) from (2) and (3) we see that to keep the number of positive discrepancies fixed under an invariant measure there should be no pair of discrepancies of different signs.

5.2. *Some auxiliary results.* Let $R_k$ be the number of distinct points visited by the random walk with transition probability $p(\cdot,\cdot)$ in $k$ steps: $R_k = \operatorname{Card}(\{X_0, \ldots, X_k\})$. Let $\tau_1 = 0$ and for $k \geq 2$ define $\tau_k = \min\{n \geq 1 : R_n = k\}$.

In [6], page 888, the following result was stated without proof.

LEMMA 5.1.
$$\sup_k \frac{\mathbb{E}(\tau_k)}{k^3} < \infty.$$

PROOF. Since $\tau_n = \sum_{k=2}^n (\tau_k - \tau_{k-1})$ for $n \geq 2$ we have
$$\frac{\mathbb{E}(\tau_n)}{n^3} \leq \sup_{2 \leq k \leq n} \mathbb{E}\left(\frac{\tau_k - \tau_{k-1}}{k^2}\right).$$

Hence, it suffices to show that

$$(33) \qquad \sup_{k \geq 2} \mathbb{E}\left(\frac{\tau_k - \tau_{k-1}}{k^2}\right) < \infty.$$

From Theorem 1 in [4] we know that there exists a constant $c$ such that for any $n$

$$\sup_x \mathbb{P}(X_n = x) \leq \frac{c}{\sqrt{n}}.$$

Thus the probability that the walk lies in a given set of $\ell$ points at time $m$ is trivially bounded above by $\ell c/\sqrt{m}$.

Let $D$ be a positive integer to be chosen later, and let $\mathcal{R}_{k-1}$ be the set of visited sites by the walk up to time $\tau_{k-1}$. We have

(34) $$\{\tau_k - \tau_{k-1} > nDk^2\} \subset \bigcap_{i=1}^n \{X_{\tau_{k-1}+iDk^2} \in \mathcal{R}_{k-1}\}.$$

But, using the strong Markov property and the previous bound, for any given $A \subset \mathbb{Z}$ with $\operatorname{Card}(A) = k-1$, we have

$$\mathbb{P}\left(\bigcap_{i=1}^n \{X_{\tau_{k-1}+iDk^2} \in A\} \Big| X_{\tau_{k-1}} \in A\right)$$

$$= \prod_{i=1}^n \mathbb{P}\left(X_{\tau_{k-1}+iDk^2} \in A \Big| \bigcap_{j=0}^{i-1} \{X_{\tau_{k-1}+jDk^2} \in A\}\right)$$

$$\leq \prod_{i=1}^n \sup_x \mathbb{P}(X_{\tau_{k-1}+iDk^2} \in A | X_{\tau_{k-1}+(i-1)Dk^2} = x)$$

$$\leq \left((k-1)\frac{c}{k\sqrt{D}}\right)^n.$$

From this bound and (34) we get

(35) $$\mathbb{P}(\tau_k - \tau_{k-1} > nDk^2) \leq \left(\frac{c}{\sqrt{D}}\right)^n.$$

Now we choose $\sqrt{D} > c$. Summing (35) over $n$ we have that $\sum_n P((\tau_k - \tau_{k-1}) > nDk^2)$ is bounded above by a constant that does not depend on $k$ which implies (33). $\square$

LEMMA 5.2. *Let $\alpha \in (0, 1/3)$. There exists a constant $c$ such that*

$$\mathbb{P}(R_k < k^\alpha) \leq \frac{ck^{3\alpha}}{k}$$

*for all $k \geq 1$.*

PROOF. For any $x \in \mathbb{R}$ denote by $\lceil x \rceil$ the smallest integer larger than or equal to $x$. Observe that

$$\{R_k < k^\alpha\} = \{\tau_{\lceil k^\alpha \rceil} > k\}.$$

By Markov's inequality

$$\mathbb{P}(\tau_{\lceil k^\alpha \rceil} > k) \leq \frac{1}{k}\mathbb{E}(\tau_{\lceil k^\alpha \rceil}).$$

From Lemma 5.1 we have that $c_0 = \sup_k \mathbb{E}(\tau_k)/k^3 < \infty$. Then

$$\mathbb{P}(R_k < k^\alpha) \leq \frac{c_0 \lceil k^\alpha \rceil^3}{k} \leq \frac{8c_0 k^{3\alpha}}{k}. \qquad \square$$

LEMMA 5.3. *For all $\rho \in (0,1)$ and $x \in \mathbb{Z}$,*

$$\sup_{z \in \mathbb{Z}} \mathbb{E}^x(\sigma(\eta_z)) \leq \mathbb{E}^x(\sigma(\eta)) \in \mathbf{L}^2(\nu_\rho)$$

*and*

$$\sup_{z \in \mathbb{Z}} \mathbb{E}^x(\sigma^*(\eta_z)) \leq \mathbb{E}^x(\sigma^*(\eta)) \in \mathbf{L}^2(\nu_\rho).$$

PROOF. Observe that for all $z \in \mathbb{Z}$, $\sigma(\eta_z) \leq \sigma(\eta)$ since it is more likely to find an empty site on $\eta_z$ than on $\eta$ (site $z$ is empty for $\eta_z$ not necessary for $\eta$). Hence, it suffices to prove that $\mathbb{E}^x(\sigma(\eta)) \in \mathbf{L}^2(\nu_\rho)$. This is done as follows:

$$\int \mathbb{E}^x(\sigma(\eta))^2 \, d\nu_\rho \leq \int \mathbb{E}^x(\sigma(\eta)^2) \, d\nu_\rho$$

$$= \int \sum_{k=1}^\infty k^2 \mathbb{P}^x(\sigma(\eta) = k) \, d\nu_\rho$$

$$\leq \int \sum_{k=1}^\infty k^2 \mathbb{E}^x\left[\prod_{\ell=1}^{k-1} \eta(X_\ell)\right] d\nu_\rho.$$

By Tonelli's theorem, this is bounded above by

$$(36) \qquad \sum_{k=1}^\infty k^2 \mathbb{E}^x(\rho^{R_{k-1}-1}) = \frac{1}{\rho}\sum_{k=0}^\infty (k+1)^2 \mathbb{E}^x(\rho^{R_k}).$$

By Lemma 5.2 there exists $c$ such that

$$\mathbb{P}(R_k < k^{1/4}) \leq ck^{-1/4}.$$

This implies:

$$\mathbb{P}(R_{\lceil k^{1/2} \rceil} < k^{1/8}) \leq \mathbb{P}(R_{\lceil k^{1/2} \rceil} < \lceil k^{1/2} \rceil^{1/4}) \leq c\lceil k^{1/2} \rceil^{-1/4} \leq ck^{-1/8}.$$

Now, for $m < n$ let $R_m^n = \text{Card}(\{X_m, \ldots, X_n\})$ and observe that, since $k > \lceil k^{1/2} \rceil(\lceil k^{1/2} \rceil - 2)$, for any $x > 0$ one has

$$\{R_k < x\} \subset \bigcap_{j=0}^{\lceil k^{1/2} \rceil - 3} \{R_{j\lceil k^{1/2} \rceil}^{(j+1)\lceil k^{1/2} \rceil} < x\}.$$

As $R_{j\lceil k^{1/2}\rceil}^{(j+1)\lceil k^{1/2}\rceil}$ for $j \geq 0$ are i.i.d. r.v.'s distributed as $R_{\lceil k^{1/2}\rceil}$, we deduce from the previous bound that

$$\mathbb{P}(R_k < k^{1/8}) < (ck^{-1/8})^{\lceil k^{1/2}\rceil - 2}.$$

Therefore conditioning on $\{R_k < k^{1/8}\}$ and its complement we obtain

$$\mathbb{E}^x(\rho^{R_k}) \leq (ck^{-1/8})^{\lceil k^{1/2}\rceil - 2} + \rho^{k^{1/8}},$$

which implies that (36) converges. $\square$

LEMMA 5.4. *For all $\rho \in (0,1)$ and all $x \in \mathbb{Z}$ we have*

$$\sum_{y \in \mathbb{Z}} |y - x|[1 - \eta(y)]\overline{q}(x, y, \eta) \in \mathbf{L}^2(\nu_\rho),$$

$$\sum_{y \in \mathbb{Z}} |y - x|[1 - \eta(y)]q(x, y, \eta) \in \mathbf{L}^2(\nu_\rho),$$

$$\sum_{y \in \mathbb{Z}} (y - x)[1 - \eta(y)]\overline{q}(x, y, \eta) = 0, \qquad \nu_\rho\text{-a.e.},$$

$$\sum_{y \in \mathbb{Z}} (y - x)[1 - \eta(y)]q(x, y, \eta) = 0, \qquad \nu_\rho\text{-a.e.}$$

PROOF. By (27)

$$\sum_{y \in \mathbb{Z}} |y - x|[1 - \eta(y)]\overline{q}(x, y, \eta) = \sum_{y \in \mathbb{Z}} |y - x|\mathbb{P}^x\left(\sum_{i=0}^{\sigma(\eta)} Y_i = y, \sigma(\eta) < \infty\right)$$

$$= \mathbb{E}^x\left(\left|\sum_{i=1}^{\sigma(\eta)} Y_i\right|, \sigma(\eta) < \infty\right).$$

By Lemma 5.3, $\sigma(\eta)$ has $\nu_\rho$-a.e. finite expectation and by the triangular inequality

$$\mathbb{E}^x\left(\left|\sum_{i=1}^{\sigma(\eta)} Y_i\right|\right) \leq \mathbb{E}^x\left(\sum_{i=1}^{\sigma(\eta)} |Y_i|\right).$$

It then follows from Wald's equation that

$$\sum_{y \in \mathbb{Z}} |y - x|[1 - \eta(y)]\overline{q}(x, y, \eta) \leq \mathbb{E}^x(\sigma(\eta)) \sum_{x \in \mathbb{Z}} |x|p(0, x), \qquad \nu_\rho\text{-a.e.}$$

The first assertion now follows from this upper bound, Lemma 5.3 and hypothesis (24). The second assertion follows from the first and the inequality $q(x,y,\eta) \leq \overline{q}(x,y,\eta)$. To prove the third assertion repeat the same argument, using hypothesis (25) to get

$$\sum_{y\in\mathbb{Z}}(y-x)[1-\eta(y)]\overline{q}(x,y,\eta) = \mathbb{E}^x\left(\sum_{i=1}^{\sigma(\eta)}Y_i, \sigma(\eta)<\infty\right)$$
$$= \mathbb{E}^x(\sigma(\eta))\sum_{x\in\mathbb{Z}}xp(0,x) = 0, \qquad \nu_\rho\text{-a.e.}$$

The last assertion is proved as the third using $\sigma(\eta_x)$ instead of $\sigma(\eta)$. $\square$

LEMMA 5.5. *Let*

$$C_\rho = \int \sum_{y\in\mathbb{Z}}|y-x|C(x,y,\eta)\,d\nu_\rho.$$

*Then, $C_\rho$ does not depend on $x$ and is finite.*

PROOF. The first assertion follows from the translation invariance of $\nu_\rho$. Hence, to complete the proof of the lemma we only need to show that

$$\int \sum_{y\in\mathbb{Z}}|y|C(0,y,\eta)\,d\nu_\rho < \infty.$$

To prove this, we first note that the integrand is bounded above by

$$\sum_{y\in\mathbb{Z}}|y|[1-\eta(y)]q(0,y,\eta) \tag{37}$$

$$+ \sum_{y\in\mathbb{Z}}|y|[1-\eta(y)]\sum_{z\neq 0, z\neq y}q(z,0,\eta)\overline{q}(0,y,\eta_z) \tag{38}$$

$$+ \sum_{y\in\mathbb{Z}}|y|q(y,0,\eta)\overline{q}(0,y,\eta_y). \tag{39}$$

By Lemma 5.4 (37) is $\nu_\rho$-integrable. We will now show that the other two terms are $\nu_\rho$-integrable. For (38) write

$$\sum_{y\in\mathbb{Z}}|y|[1-\eta(y)]\sum_{z\neq 0, z\neq y}q(z,0,\eta)\overline{q}(0,y,\eta_z)$$
$$= \sum_{z\neq 0}q(z,0,\eta)\sum_{y\neq z}|y|[1-\eta(y)]\overline{q}(0,y,\eta_z).$$

Since for $y\neq z$, $\overline{q}(0,y,\eta_z) \leq \overline{q}(0,y,\eta)$, this last expression is bounded above by

$$\left(\sum_{z\neq 0}q(z,0,\eta)\right)\left(\sum_{y\in\mathbb{Z}}|y|[1-\eta(y)]\overline{q}(0,y,\eta)\right).$$

Hence, by Lemma 5.4 it is now enough to show that $\sum_{z \neq 0} q(z, 0, \eta)$ is in $\mathbf{L}^2(\nu_\rho)$. For this purpose, note that $q(z, 0, \eta) = q^*(0, z, \eta)$ and $\sum_{z \neq 0} q^*(0, z, \eta)$ is bounded above by $\mathbb{E}(\sigma^*(\eta_0))$ which is in $\mathbf{L}^2(\nu_\rho)$ by Lemma 5.3.

For (39) first note that since $q(y, 0, \eta) \leq 1$, it suffices to show that $\sum_{y \in \mathbb{Z}} |y| \bar{q}(0, y, \eta_y) = \sum_{y \in \mathbb{Z}} |y| \bar{q}(0, y, \eta)$ is $\nu_\rho$-integrable. Since $\nu_\rho$ is a product measure and $\bar{q}(0, y, \eta)$ does not depend on $\eta(y)$, we have

$$\int \sum_{y \in \mathbb{Z}} |y| [1 - \eta(y)] \bar{q}(0, y, \eta) \, d\nu_\rho = (1 - \rho) \int \sum_{y \in \mathbb{Z}} |y| \bar{q}(0, y, \eta) \, d\nu_\rho.$$

Since $\rho < 1$ and the left-hand side above is finite by Lemma 5.4 we also have

$$\int \sum_{y \in \mathbb{Z}} |y| \bar{q}(0, y, \eta) \, d\nu_\rho < \infty$$

which completes the proof. $\square$

LEMMA 5.6. *For all $\rho \in (0, 1)$ and all $x \in \mathbb{Z}$*

$$\sum_{y \in \mathbb{Z}} (y - x) B(x, y, \eta, \xi) = 0, \qquad \widetilde{\nu}_\rho - a.e.$$

PROOF. Since $B(x, y, \eta, \xi) \leq C(x, y, \eta)$ the series converges absolutely $\widetilde{\nu}_\rho$-a.e. by Lemma 5.5. Therefore it suffices to show that

$$\sum_{y \in \mathbb{Z}} (y - x) \mathbf{I}_{\eta(x)=1, \xi(x)=\eta(y)=0} \, q(x, y, \eta) = 0, \qquad \widetilde{\nu}_\rho\text{-a.e.}$$

and that for all $z \neq x$

$$q(z, x, \eta\xi) \mathbf{I}_{\eta(z)=\xi(z)=1} \sum_{y \in \mathbb{Z}} (y - x) \mathbf{I}_{\eta(x)=1, \xi(x)=\eta_z(y)=0} \, \bar{q}(x, y, \eta_z) = 0, \qquad \widetilde{\nu}_\rho\text{-a.e.}$$

But the first of these terms is equal to

$$\mathbf{I}_{\eta(x)=1, \xi(x)=0} \sum_{y \in \mathbb{Z}} (y - x)[1 - \eta(y)] q(x, y, \eta),$$

and the second is equal to

$$q(z, x, \eta\xi) \mathbf{I}_{\eta(z)=\xi(z)=\eta(x)=1, \xi(x)=0} \sum_{y \in \mathbb{Z}} (y - x)[1 - \eta_z(y)] \bar{q}(x, y, \eta_z),$$

and both these expressions vanish $\nu_\rho$-a.e. by Lemma 5.4 and the fact that if $\eta$ is distributed according to $\nu_\rho$, then the distribution of $\eta_z$ is absolutely continuous with respect to $\nu_\rho$. Since this measure is the first marginal of $\widetilde{\nu}_\rho$ the lemma follows. $\square$

5.3. *Proof of Theorem* 1.5. Let $\widetilde{\nu}_\rho$ be an invariant measure for the coupled process with respective marginals $\nu_\rho$, $\rho \in (0,1)$ and $\mu \in \mathcal{I}$.

Let $f_n(\eta, \xi) = \sum_{x=-n}^{n} [\eta(x) - \xi(x)]^+$. This means that $f_n$ counts the number of positive discrepancies in the interval $[-n, n]$.

LEMMA 5.7. *The sum of the positive terms of* $\widetilde{L}f_n(\eta, \xi)$ *is equal to*

$$\sum_{\substack{|x|>n \\ |y|\leq n}} A(x,y,\eta,\xi),$$

*and the sum of the negative terms is less than or equal to*

$$-\sum_{\substack{|x|\leq n \\ |y|>n}} B(x,y,\eta,\xi) - \sum_{\substack{-n\leq x,y\leq n \\ x\neq y}} D(x,y,\eta,\xi).$$

PROOF. We compute $\widetilde{L}f_n(\eta, \xi)$ using the expression given in Section 1. Note that in this case (2), (5), (9), (12) and (13) vanish and that the positive terms only come from (4), (6) and (7).

The sum of the positive terms of (4) is equal to

$$\sum_{\substack{x\in\mathbb{Z}, |y|>n \\ |z|\leq n}} \mathbf{I}_{\eta(x)=\xi(x)=\eta(y)=1, \xi(y)=\eta(z)=\xi(z)=0}\, q(x,y,z,\xi,\eta),$$

which using the notation of Section 5.1 we can write as

(40) $$\sum_{\substack{|y|>n \\ |z|\leq n}} \mathbf{I}_{\eta(y)=1,\xi(y)=0} \sum_{x\neq y,z} \mathbf{I}_{\eta(x)=\xi(x)=1,\eta_x(z)=\xi_x(z)=0}\, q(x,y,\xi\eta)\overline{q}(y,z,\eta_x).$$

The positive terms from (6) are equal to

$$\sum_{\substack{|x|\leq n \\ |y|>n}} \mathbf{I}_{\eta(x)=\xi(x)=\eta(y)=1, \xi(y)=0}\, q(x,y,\xi,\eta)$$

(41)
$$= \sum_{\substack{|x|\leq n \\ |y|>n}} \mathbf{I}_{\eta(y)=1,\xi(y)=0}\mathbf{I}_{\eta(x)=\xi(x)=1,\eta_x(x)=\xi_x(x)=0}\, q(x,y,\xi\eta)\overline{q}(y,x,\eta_x).$$

Finally from (7) we get

$$\sum_{\substack{|x|>n \\ |y|\leq n}} \mathbf{I}_{\eta(x)=1,\xi(x)=\eta(y)=\xi(y)=0}\, q(x,y,\eta)$$

(42)
$$= \sum_{\substack{|x|>n \\ |y|\leq n}} \mathbf{I}_{\eta(x)=1,\xi(x)=0}\mathbf{I}_{\eta(y)=\xi(y)=0}\, q(x,y,\eta).$$

Applying the permutation $y \to x \to z \to y$ to (40) and the permutation $y \to x \to y$ to (41) and then adding (40), (41) and (42) we get

$$\sum_{\substack{|x|>n \\ |y|\leq n}} A(x,y,\eta,\xi),$$

which proves the first assertion of the lemma.

Among the negative terms we focus on (part of) those obtained through (4), (6) and (7).

Part of the negative contribution of (4) is given by

$$\sum_{\substack{x\in\mathbb{Z},|y|\leq n \\ |z|>n}} \mathbf{I}_{\eta(x)=\xi(x)=\eta(y)=1,\xi(y)=\eta(z)=0}\, q(x,y,z,\xi,\eta),$$

which using the notation of Section 5.1 we can write as

(43) $$\sum_{\substack{|y|\leq n \\ |z|>n}} \mathbf{I}_{\eta(y)=1,\xi(y)=0} \sum_{x\neq y,z} \mathbf{I}_{\eta(x)=\xi(x)=1,\eta_x(z)=0}\, q(x,y,\xi\eta)\overline{q}(y,z,\eta_x).$$

Observe that this is just part of the negative contribution of (4) since there should be another term with $x,y$ chosen as before and $z \in [-n,n]$ with $\xi(z) = 1$.

From (6) we get as negative terms

(44)
$$\sum_{\substack{|x|>n \\ |y|\leq n}} \mathbf{I}_{\eta(x)=\xi(x)=\eta(y)=1,\xi(y)=0}\, q(x,y,\xi,\eta)$$
$$= \sum_{\substack{|x|>n \\ |y|\leq n}} \mathbf{I}_{\eta(y)=1,\xi(y)=0}\mathbf{I}_{\eta(x)=\xi(x)=1,\eta_x(x)=0}\, q(x,y,\xi\eta)\overline{q}(y,x,\eta_x).$$

Finally, from (7) we separate the contribution in two sums

$$\sum_{\substack{|x|\leq n \\ |y|>n}} \mathbf{I}_{\eta(x)=1,\xi(x)=\eta(y)=0}\, q(x,y,\eta)$$

(45) $$+ \sum_{\substack{|x|\leq n \\ |y|\leq n}} \mathbf{I}_{\eta(x)=\xi(y)=1,\xi(x)=\eta(y)=0}\, q(x,y,\eta)$$

$$= \sum_{\substack{|x|\leq n \\ |y|>n}} \mathbf{I}_{\eta(x)=1,\xi(x)=0}\mathbf{I}_{\eta(y)=0}\, q(x,y,\eta) + \sum_{\substack{|x|\leq n \\ |y|\leq n, x\neq y}} D(x,y,\eta,\xi).$$

Now, relabelling when necessary and summing the expressions in (43), (44) and the first sum in (45) we get

$$\sum_{\substack{|x|\leq n \\ |y|>n}} B(x,y,\eta,\xi),$$

which proves the second assertion of the lemma. □

It now follows from Theorem 1.2 and Lemma 5.7 that

$$\sum_{\substack{x<-n \\ y\geq -n}} A_\rho(x,y) + \sum_{\substack{x>n \\ y\leq n}} A_\rho(x,y) - \sum_{\substack{|x|\leq n \\ |y|>n}} B_\rho(x,y) \geq \sum_{\substack{-n\leq x,y\leq n \\ x\neq y}} D_\rho(x,y),$$

which by (32) implies

$$\sum_{\substack{x<-n \\ y\geq -n}} B_\rho(x,y) + \sum_{\substack{x>n \\ y\leq n}} B_\rho(x,y) - \sum_{\substack{|x|\leq n \\ |y|>n}} B_\rho(x,y) \geq \sum_{\substack{-n\leq x,y\leq n \\ x\neq y}} D_\rho(x,y).$$

Call $\Gamma_\rho(n)$ the left-hand side of the previous inequality. We will show that $\limsup_N \frac{1}{N+1}\sum_{n=0}^N \Gamma_\rho(n) \leq 0$. This will imply that

$$(46) \qquad D_\rho(x,y) = 0 \qquad \text{for all } x,y.$$

Let $M = M(N)$ be such that $0 < M < N$, $M \to \infty$ as $N \to \infty$ and

$$(47) \qquad \lim_{N\to\infty} \frac{M}{N} = 0.$$

Then

$$\frac{1}{N+1}\sum_{n=0}^N \Gamma_\rho(n)$$

$$(48) \qquad \leq \frac{1}{N+1}\Bigg[\sum_{x=-N+M+1}^{-1}\sum_{y=x+1}^{\infty}(y-x)B_\rho(x,y)$$

$$(49) \qquad + \sum_{x=-N-M}^{-N+M}\sum_{y=x+1}^{\infty}(y-x)B_\rho(x,y)$$

$$(50) \qquad + \sum_{x=-\infty}^{-N-M-1}\sum_{y=-N}^{\infty}[(y-x)\wedge(N+1)]B_\rho(x,y)$$

$$(51) \qquad - \sum_{x=-N}^{-1}\sum_{y=-\infty}^{x-1}[(x-y)\wedge(N+x+1)]B_\rho(x,y)\Bigg]$$

$$(52) \quad + \frac{1}{N+1}\Bigg[\sum_{x=1}^{N-M-1}\sum_{y=-\infty}^{x-1}(x-y)B_\rho(x,y)$$

$$(53) \quad + \sum_{x=N-M}^{N+M}\sum_{y=-\infty}^{x-1}(x-y)B_\rho(x,y)$$

$$(54) \quad + \sum_{x=N+M+1}^{\infty}\sum_{y=-\infty}^{N}[(x-y)\wedge(N+1)]B_\rho(x,y)$$

$$(55) \quad - \sum_{x=1}^{N}\sum_{y=x+1}^{\infty}[(y-x)\wedge(N-x+1)]B_\rho(x,y)\Bigg].$$

Symbolically we will write the right-hand side as

$$\frac{1}{N+1}[(48)+(49)+(50)-(51)] + \frac{1}{N+1}[(52)+(53)+(54)-(55)].$$

The next three lemmas will show that

$$\limsup_N \frac{1}{N+1}[(48)+(49)+(50)-(51)] \leq 0.$$

The second bracket is treated analogously.

LEMMA 5.8. $\frac{1}{N+1}(49) \to 0$ as $N \to \infty$.

PROOF. It follows from Lemma 5.5 that

$$(49) \leq \sum_{x=-N-M}^{-N+M}\sum_{y\in\mathbb{Z}}|y-x|C_\rho(x,y) \leq (2M+1)C_\rho.$$

Hence, the lemma is a consequence of (47). □

LEMMA 5.9. $\frac{1}{N+1}(50) \to 0$ as $N \to \infty$.

PROOF. Observe that

$$\frac{1}{N+1}(50) \leq \sum_{x=-\infty}^{-N-M-1}\sum_{y=-N}^{\infty}B_\rho(x,y)$$

$$\leq \sum_{x=-\infty}^{-N-M-1}\sum_{y=-N}^{\infty}C_\rho(x,y)$$

$$\leq \sum_{k=M+1}^{\infty}kC_\rho(0,k),$$

which goes to 0 as $M \to \infty$ since by Lemma 5.5 this series converges. □

LEMMA 5.10. $\frac{1}{N+1}[(48) - (51)] \to 0$ as $N \to \infty$.

PROOF. The expression (48)–(51) is bounded above by

$$\sum_{x=-N+M+1}^{-1} \sum_{y=-N+1}^{\infty} (y-x)B_\rho(x,y).$$

For $x \in (-N+M, 0)$, by Lemma 5.6 we have

$$\sum_{y=-N+1}^{\infty} (y-x)B_\rho(x,y) = \sum_{y \in \mathbb{Z}} (y-x)B_\rho(x,y) - \sum_{y=-\infty}^{-N} (y-x)B_\rho(x,y)$$

$$= 0 + \sum_{y=-\infty}^{-N} (x-y)B_\rho(x,y).$$

Hence, for these $x$'s it follows from (32) that

$$\sum_{y=-N+1}^{\infty} (y-x)B_\rho(x,y) \leq \sum_{y=-\infty}^{-N} (x-y)C_\rho(x,y) \leq \sum_{k=M}^{\infty} kC_\rho(k,0).$$

Therefore

$$\frac{1}{N+1}[(48) - (51)] \leq \frac{N-M}{N+1} \sum_{k=M}^{\infty} kC_\rho(k,0) \leq \sum_{k=M}^{\infty} kC_\rho(0,-k),$$

which goes to 0 as $M \to \infty$ by Lemma 5.5. This concludes the proof of Lemma 5.10. □

These lemmas imply (46) which implies (23). As explained at the beginning of this section this proves Theorem 1.5.

Suppose for the final part of this section that $\mathbf{X} = \{0,1\}^{\mathbb{Z}^d}$ and $p$ is an irreducible random walk on $\mathbb{Z}^d$.

SKETCH OF PROOF OF THEOREM 1.4. Suppose $\mu \in (\mathcal{I} \cap \mathcal{S})$ puts no mass on **1**. Apply Theorem 1.2 to the function $f(\eta, \xi) = [\eta(x) - \xi(x)]^+$ and the translation-invariant coupled measure $\widetilde{\nu}_\rho$ of Lemma 4.1 with marginals $\nu_\rho$ and $\mu$. Due to the translation invariance of $p(x,y)$ many of the terms of $\int \widetilde{L}f \, d\widetilde{\nu}_\rho$ cancel each other and only negative terms are left. Since the sum of these is equal to $-D_\rho$ we again obtain (23) for any $\rho \in (0,1)$. □

**6. The nearest-neighbor case.** In this section we suppose that $d = 1$, $p(x, x+1) = p$ and $p(x, x-1) = q$, $p + q = 1$, and prove Theorem 1.6. In view of Theorem 1.5 we may assume $p \neq 1/2$ and by symmetry it suffices to prove the theorem for $p > 1/2$. Initially, we assume that $p < 1$ to keep the transition matrix irreducible. At the end of the section we will explain which modifications are needed to treat the case $p = 1$. As in the previous section let

$$f_n(\eta, \xi) = \sum_{x=-n}^{n} [\eta(x) - \xi(x)]^+$$

be the number of positive discrepancies between $\eta$ and $\xi$ in $[-n, n]$ and let

$$g_n(\eta, \xi) = \text{Card}(\{x \in [-n, n) : \eta(x) - \xi(x) = 1 \text{ and } \exists y \in (x, n]$$
$$\text{such that } \forall z \in (x, y), \eta(z) - \xi(z) = 0$$
$$\text{and } \eta(y) - \xi(y) = -1\}).$$

This means that $g_n$ counts the number of times we see $\eta - \xi$ changing from $+1$ to $-1$ when we move from $-n$ to $n$.

Let $\mu \in \mathcal{I}$ be such that $\mu(\{\mathbf{1}\}) = 0$, let $\rho \in (0, 1)$ and let $\widetilde{\nu}_\rho$ be a probability measure on $\mathbf{X} \times \mathbf{X}$ with marginals $\nu_\rho$ and $\mu$ which is invariant for the coupled process. As explained at the end of Section 4, to prove Theorem 1.6 it suffices to show that $\widetilde{\nu}_\rho$ is such that

$$\widetilde{\nu}_\rho(\{(\eta, \xi) : \eta(x) = \xi(y) = 1, \eta(y) = \xi(x) = 0\}) = 0,$$

whenever $p(x, y) > 0$. The proof relies loosely on the following observations:

(1) Suppose $g_n > 1$; then it can decrease with time.
(2) To compensate this, some boundary effects must increase $g_n$.
(3) These boundary effects require $g_n$ to grow linearly with $n$.

Hence, many discrepancies of opposite signs are not far from each other and discrepancies disappear in the interval $[-n, n]$ at a rate which grows linearly with $n$. But positive discrepancies enter that interval at a bounded rate and we get a contradiction.

At different parts of the proof we must follow the movement of discrepancies. When doing so it is important to realize that a discrepancy at a site $x$ can move if the clock rings at that site, but it can also move if the clock rings at another site $z$ occupied by both $\eta$ and $\xi$ particles because these two particles may arrive at $x$ where only one will stay. Since we are dealing with a nearest-neighbor random walk, for this last case to occur, all the sites between $x$ and $z$ must be occupied by both $\eta$ and $\xi$ particles.

The generator $\widetilde{L}$ of the coupled process will be applied to several cylinder functions $f$. We denote by $\widetilde{L}^+ f$ ( $-\widetilde{L}^- f$) the sum of the positive (negative) terms of $\widetilde{L}f$. Hence $\widetilde{L}f = \widetilde{L}^+ f - \widetilde{L}^- f$.

PROPOSITION 6.1. *Let $\widetilde{\nu}_\rho$ be as above. Then for any $n \geq 1$*

$$\widetilde{\nu}_\rho\{(\eta, \xi) : g_n(\eta, \xi) \geq 2\} = 0.$$

To prove the proposition let $R(x, \widetilde{\nu}_\rho)$ be the average under $\widetilde{\nu}_\rho$ of the rate at which a positive discrepancy originally at $x$ disappears either because it moves to a site occupied by a negative discrepancy or because a negative discrepancy moves to $x$.

LEMMA 6.1. *Fix $0 < \rho < 1$, let $\ell \in \mathbb{N}$ and let $\widetilde{\nu}_\rho$ be as in Proposition 6.1. Then, there exists an increasing function $h_\ell : [0,1) \to [0, \infty)$ such that:*

(i) $h_\ell(t) > 0$ *if* $t > 0$, *and*
(ii) *for any $x \in \mathbb{Z}$,*

$$\sum_{y=x-\ell}^{x} R(y, \widetilde{\nu}_\rho)$$
$$\geq h_\ell(\widetilde{\nu}_\rho\{(\eta, \xi) : \eta(x-\ell) - \xi(x-\ell) = 1, \eta(x-i) - \xi(x-i) = 0,$$
$$1 \leq i \leq \ell - 1, \eta(x) - \xi(x) = -1\}).$$

PROOF. We prove the lemma by induction on $\ell$. For $\ell = 1$ the result is trivial with $h_1(c) = c$. We now assume that the conclusion of the lemma holds for $\ell = 1, \ldots, k-1$, $k \geq 2$, and suppose that

$$\widetilde{\nu}_\rho\{(\eta, \xi) : \eta(x-k) - \xi(x-k) = 1, \eta(x-i) - \xi(x-i) = 0,$$
$$1 \leq i \leq k-1, \eta(x) - \xi(x) = -1\} = c > 0. \tag{56}$$

Let $\overline{p}$ be the probability that a $p, q$ random walk starting at $0$ never returns to the origin (and therefore drifts to $+\infty$). It follows from (56) that either

$$\widetilde{\nu}_\rho\{(\eta, \xi) : \eta(x-k) - \xi(x-k) = 1, \eta(x-i) = \xi(x-i) = 1,$$
$$1 \leq i \leq k-1, \eta(x) - \xi(x) = -1\} \geq \frac{c}{k},$$

or $\exists r$, $1 \leq r \leq k-1$ such that

$$\widetilde{\nu}_\rho\{(\eta, \xi) : \eta(x-k) - \xi(x-k) = 1, \eta(x-i) = \xi(x-i) = 1,$$
$$r+1 \leq i \leq k-1, \eta(x-r) = \xi(x-r) = 0, \tag{57}$$
$$\eta(x-j) = \xi(x-j), 1 \leq j \leq r-1, \eta(x) - \xi(x) = -1\} \geq \frac{c}{k}.$$

In the first case the average under $\widetilde{\nu}_\rho$ of the rate at which a positive discrepancy at $x - k$ jumps to $x$ (and is annihilated) is at least $\overline{p} \cdot c/k$.

To treat the second case, let $r \in [1, k-1]$ be such that (57) is satisfied. We apply Theorem 1.2 to the function

$$t_r(\eta, \xi) = \mathbf{I}_{\{\eta(x-r) - \xi(x-r) = 1, \eta(x-j) = \xi(x-j), 1 \leq j \leq r-1, \eta(x) - \xi(x) = -1\}}.$$

Among the positive terms of $\int \widetilde{L} t_r(\eta, \xi) \, d\widetilde{\nu}_\rho$ there is one due to the jump of a positive discrepancy from $x - k$ to $x - r$. This term is greater than or equal to $\overline{p} \cdot c/k$. Hence, by Theorem 1.2, $\int \widetilde{L}^- t_r(\eta, \xi) \, d\widetilde{\nu}_\rho \geq \overline{p} \cdot c/k$. We now separate the terms of $\widetilde{L}^- t_r$ in five groups:

(1) terms due to jumps of particles from sites $y > x$ into the interval $[x - r, x]$;

(2) terms due to jumps starting from the interval $[x - r, x]$;

(3) terms due to jumps of an $\eta$ particle (and possibly a $\xi$ particle too) from a site $y \in [x - r - n_0, x - r)$ into the interval $[x - r, x]$ ($n_0$ will be chosen later);

(4) terms due to jumps of an $\eta$ particle (and possibly a $\xi$ particle too) from a site $y < x - r - n_0$ into the interval $[x - r, x]$;

(5) terms due to jumps of a $\xi$ particle only from a site $y < x - r$ into the interval $[x - r, x]$ that annihilate the positive discrepancy at $x - r$.

Let $\overline{q} = q/p$, then $(\overline{q})^n$ is the probability that a $p, q$ random walk starting at $n$ ever hits 0. The sum of the terms in group 1 is bounded above by

$$\sum_{n=1}^{\infty} (\overline{q})^n \int t_r(\eta, \xi) \, d\widetilde{\nu}_\rho = \frac{\overline{q}}{1 - \overline{q}} \int t_r(\eta, \xi) \, d\widetilde{\nu}_\rho.$$

The sum of the terms in group 2 is bounded above by

$$(r + 1) \int t_r(\eta, \xi) \, d\widetilde{\nu}_\rho.$$

The sum of the terms in group 3 is bounded above by

$$n_0 \int t_r(\eta, \xi) \, d\widetilde{\nu}_\rho.$$

The sum of the terms in group 4 is bounded above by

$$\sum_{n=n_0}^{\infty} \rho^n.$$

We choose $n_0$ in such a way that this sum is $< \overline{p} c/5(k)$. Then, either the jumps of group 5 occur at a rate greater than or equal to $\overline{p} c/(5k)$, or one of the first three upper bounds is at least $\overline{p} c/(5k)$. In this second case we get

$$\int t_r(\eta, \xi) \, d\widetilde{\nu}_\rho \geq \frac{c\overline{p}}{5k} \min\left\{ \frac{1 - \overline{q}}{\overline{q}}, \frac{1}{r + 1}, \frac{1}{n_0} \right\} := A_c.$$

Therefore taking

$$h_k(c) = \min\left\{ \frac{c\overline{p}}{5k}, h_1(A_c), h_2(A_c), \ldots, h_{k-1}(A_c) \right\},$$

the result follows from the inductive hypothesis. $\square$

LEMMA 6.2. *Suppose $\widetilde{\nu}_\rho$ is as in Proposition 6.1. If for some $n_0 \in \mathbb{N}$*

$$\widetilde{\nu}_\rho(\{(\eta,\xi) : g_{n_0}(\eta,\xi) \geq 2\}) = \delta_1 > 0,$$

*then there exists a $\delta_2 > 0$ such that $\forall n \geq n_0$*

$$\int \widetilde{L}^- g_n \, d\widetilde{\nu}_\rho \geq \delta_2.$$

The proof of this lemma requires some extra notation: For fixed $(\eta,\xi)$ such that $g_{n_0}(\eta,\xi) = m \geq 2$ let

$$r_1(\eta,\xi) < r_2(\eta,\xi) < \cdots < r_{2m+1}(\eta,\xi)$$

be the smallest integers in the interval $[-n_0, n_0+1]$ such that:

(i) each of the intervals $[r_1, r_2-1], \ldots, [r_{2m-1}, r_{2m}-1]$ contains at least one positive discrepancy and no negative discrepancies,
(ii) each of the intervals $[r_2, r_3-1], \ldots, [r_{2m}, r_{2m+1}-1]$ contains at least one negative discrepancy and no positive discrepancies.

For each $i \in \{1, \ldots, m-1\}$, let

$$x_{i,1} > \cdots > x_{i,j_i}$$

be the position of the negative discrepancies in the interval $[r_{2i}, r_{2i+1}-1]$ and let

$$y_{i,1} < \cdots < y_{i,j'_i}$$

be the position of the positive discrepancies in the interval $[r_{2i+1}, r_{2i+2}-1]$. Define

$$N_i(\eta,\xi) = \sum_{r=1}^{j_i \wedge j'_i} (y_{i,r} - x_{i,r}) \quad \text{and} \quad N(\eta,\xi) = \min_{1 \leq i \leq m-1} N_i(\eta,\xi).$$

The quantity $N_i(\eta,\xi)$ is the total distance that the discrepancies in the interval $[r_{2i}, r_{2i+1}-1]$ must travel to the right to cause $g_n$ ($n \geq n_0$) to decrease.

PROOF OF LEMMA 6.2. Let

$$A = \{(\eta,\xi) : g_{n_0}(\eta,\xi) \geq 2\}.$$

For $k \geq 1$, let

$$A_k = \{(\eta,\xi) \in A : N(\eta,\xi) = k\}$$

and

$$B_k = \bigcup_{i=1}^k A_i.$$

Since $N(\eta, \xi) \leq 4n_0^2$, it follows from the hypothesis of the lemma that for some $k \leq 4n_0^2$ we have

$$\widetilde{\nu}_\rho(A_k) \geq \frac{\delta_1}{4n_0^2}.$$

Hence, to prove the lemma it suffices to show that for each $k \geq 1$, there exists $c_k > 0$ such that

(58) $$\int \widetilde{L}^- g_n(\eta, \xi) \, d\widetilde{\nu}_\rho \geq c_k \widetilde{\nu}_\rho(A_k) \qquad \forall n \geq n_0.$$

We proceed by induction on $k$. On the set $A_1$ there is a negative discrepancy in some site $x$ in the interval $[-n_0, n_0]$ whose jump to the right causes $g_n$ to decrease. Hence, for any configuration in $A_1$, $g_n$ decreases at rate at least $p$ and (58) holds for $k = 1$ and $c_1 = p$. Suppose now that it holds for $k = 1, \ldots, k_0 - 1$. Let $(\eta, \xi) \in A_{k_0}$ and let $\bar{i}$ be such that $N_{\bar{i}}(\eta, \xi) = N(\eta, \xi)$. Then, at rate at least $\bar{p}$ (as defined in the proof of Lemma 6.1) the negative discrepancy at $x_{\bar{i},1}$ moves to the right. When this happens either $g_n$ decreases or the new configuration belongs to $B_{k_0-1}$. Therefore, for all $n \geq n_0$ either

$$\int \widetilde{L}^+ \mathbf{I}_{B_{k_0-1}} \, d\widetilde{\nu}_\rho \geq \frac{\bar{p}}{2} \widetilde{\nu}_\rho(A_{k_0})$$

or

$$\int \widetilde{L}^- g_n \, d\widetilde{\nu}_\rho \geq \frac{\bar{p}}{2} \widetilde{\nu}_\rho(A_{k_0}).$$

In the second case the proof is complete. In the first case, by Theorem 1.2

$$\int \widetilde{L}^- \mathbf{I}_{B_{k_0-1}} \, d\widetilde{\nu}_\rho \geq \frac{\bar{p}}{2} \widetilde{\nu}_\rho(A_{k_0}).$$

Write $\widetilde{L}^- \mathbf{I}_{B_{k_0-1}}$ as $\widetilde{L}_1^- \mathbf{I}_{B_{k_0-1}} + \widetilde{L}_2^- \mathbf{I}_{B_{k_0-1}} + \widetilde{L}_3^- \mathbf{I}_{B_{k_0-1}}$, where $L_1^-$, $L_2^-$ and $L_3^-$ contain the terms of $\widetilde{L}^-$ due to jumps of particles at sites in $(-\infty, -n_0 - 1]$, $[-n_0, n_0]$ and $[n_0 + 1, \infty)$, respectively. It follows from the last inequality that either

$$\int \widetilde{L}_2^- \mathbf{I}_{B_{k_0-1}} \, d\widetilde{\nu}_\rho + \int \widetilde{L}_3^- \mathbf{I}_{B_{k_0-1}} \, d\widetilde{\nu}_\rho \geq \frac{\bar{p}}{4} \widetilde{\nu}_\rho(A_{k_0}),$$

or

$$\int \widetilde{L}_1^- \mathbf{I}_{B_{k_0-1}} \, d\widetilde{\nu}_\rho \geq \frac{\bar{p}}{4} \widetilde{\nu}_\rho(A_{k_0}).$$

In the first case we use the inequalities

$$\int \widetilde{L}_2^- \mathbf{I}_{B_{k_0-1}} \, d\widetilde{\nu}_\rho \leq (2n_0 + 1) \widetilde{\nu}_\rho(B_{k_0-1})$$

and

$$\int \tilde{L}_3^- \mathbf{I}_{B_{k_0-1}} \, d\tilde{\nu}_\rho \leq \sum_{n=1}^{\infty} \overline{q}^n \tilde{\nu}_\rho(B_{k_0-1}) = \overline{q}(1-\overline{q})^{-1} \tilde{\nu}_\rho(B_{k_0-1})$$

to obtain

$$\tilde{\nu}_\rho(B_{k_0-1}) \geq \frac{\overline{p}}{4} \frac{1}{2n_0 + 1 + \overline{q}(1-\overline{q})^{-1}} \tilde{\nu}_\rho(A_{k_0}),$$

and the result follows from the inductive hypothesis. To complete the proof in the second case it suffices to check all the terms of $L_1^- \mathbf{I}_{B_{k_0-1}}$ are also terms of $L^- g_n$ for all $n \geq n_0$. This is done as follows: The arrival of an $\eta$ particle from the left of $-n_0$ can increase the value of $N(\eta, \xi)$ only if it destroys a full interval of negative discrepancies situated to the right of the first positive discrepancy and the arrival of $\xi$ particle can increase the value of $N(\eta, \xi)$ only if it destroys the first interval of positive discrepancies. In both cases the value of $g_n$ decreases for all $n \geq n_0$. □

The following lemma is elementary:

LEMMA 6.3. *Let $\rho$ be a strictly positive real number, let $n, m \in \mathbb{N}$ and let $a_i$, $i = 1, \ldots, m$, be real numbers such that:*

(i) $0 \leq a_i \leq 1$, $\forall 1 \leq i \leq m$,
(ii) $\sum_{i=1}^m a_i \geq \rho n$.

*Then,*

$$\mathrm{Card}\left(\left\{i : 1 \leq i \leq m, \ a_i \geq \frac{\rho}{2m} n\right\}\right) \geq \frac{\rho}{2} n.$$

PROOF. Write

$$\rho n \leq \sum_{i=1}^m a_i = \sum_{1 \leq i \leq m : a_i < (\rho/(2m))n} a_i + \sum_{1 \leq i \leq m : a_i \geq (\rho/(2m))n} a_i.$$

Since the first sum in the right-hand side is bounded above by $\frac{\rho}{2} n$, the second sum is bounded below by $\frac{\rho}{2} n$ and the result follows because all the $a_i$'s are in $[0, 1]$. □

Before stating our last lemma we introduce some notation: For $k, n \in \mathbb{N}$ we say that an increasing sequence of integers $x_1, y_1, \ldots, x_m, y_m$ is the $k$-partition of the interval $[-n, n]$ if:

(i) $x_1 = -n$, $y_m = n$,
(ii) $x_{i+1} = y_i + 1$ for $i = 1, \ldots, m-1$,

(iii) $y_i - x_i = k$ for $i = 1, \ldots, m-1$,
(iv) $0 \leq y_m - x_m \leq k$.

For integers $a < b$ let $g_{a,b}(\eta, \xi)$ be as $g_n(\eta, \xi)$ with $[a,b]$ substituting for $[-n,n]$. It follows from these definitions that if $x_1, y_1, \ldots, x_m, y_m$ is the $k$-partition of the interval $[-n,n]$, then

$$(59) \qquad g_n(\eta, \xi) \leq \sum_{i=1}^m g_{x_i, y_i}(\eta, \xi) + m - 1.$$

LEMMA 6.4. *Let $\nu$ be a probability measure on $\mathbf{X} \times \mathbf{X}$ such that for some $n_0 \in \mathbb{N}$ and $0 < \alpha \leq 1$*

$$\int g_n(\eta, \xi) \, d\nu \geq \alpha n \qquad \forall n \geq n_0.$$

*Let $k = [4/\alpha] + 1$ and let $x_1, y_1, \ldots, x_m, y_m$ be the $k$-partition of $[-n, n]$. Then, there exists $\beta > 0$ such that*

$$\mathrm{Card}(\{1 \leq i \leq m-1 : \nu(g_{x_i, y_i}(\eta, \xi) \geq 1) \geq \beta\}) \geq \beta n \qquad \forall n \geq n_0.$$

PROOF. Since $(m-1)k \leq 2n$ and $k > 4/\alpha$ we must have $m - 1 \leq 2^{-1} \alpha n$. It then follows from (59) that

$$\sum_{i=1}^m \int g_{x_i, y_i}(\eta, \xi) \, d\nu \geq \frac{\alpha}{2} n \qquad \forall n \geq n_0.$$

Since $g_{x_i, y_i}$ is bounded above by $k \leq 5/\alpha$, this implies

$$\sum_{i=1}^m \int g_{x_i, y_i}(\eta, \xi) \wedge 1 \, d\nu \geq \frac{\alpha}{5} \sum_{i=1}^m \int g_{x_i, y_i}(\eta, \xi) \, d\nu \geq \frac{\alpha^2}{10} n.$$

Since the last left-hand side can be written as $\sum_{i=1}^m \nu(g_{x_i, y_i} \geq 1)$ the lemma follows from the previous lemma and the fact that $n/m \geq 1$. □

PROOF OF PROPOSITION 6.1. We argue by contradiction: If the proposition does not hold, then there exists a $\delta_1 > 0$ and $n_0 \in \mathbb{N}$ such that

$$\widetilde{\nu}_\rho \{(\eta, \xi) : g_n(\eta, \xi) \geq 2\} \geq \delta_1 \qquad \forall n \geq n_0.$$

By Lemma 6.2 there exists a $\delta_2 > 0$ such that for all $n \geq n_0$, $\int \widetilde{L}^- g_n \, d\widetilde{\nu}_\rho \geq \delta_2$. Therefore, by Theorem 1.2 $\int \widetilde{L}^+ g_n \, d\widetilde{\nu}_\rho \geq \delta_2$. The terms of $\widetilde{L}^+ g_n$ are due either to positive discrepancies coming from the left of $-n$ or to negative discrepancies coming from the right of $n$.

A positive discrepancy at $-n - k$ may move into the interval $[-n, n]$ if all the sites between $-n - k$ and $-n$ are occupied by $\eta$ particles and either because the clock rings at $-n - k$ or the clock rings at another site $x$ occupied by both $\eta$ and $\xi$ particles and all the sites between $x$ and $-n-k$ are also occupied by both $\eta$ and $\xi$ particles. Since the $\eta$ particles are distributed according to $\nu_\rho$, the contribution of this discrepancy to $\int \widetilde{L}^+ g_n \, d\widetilde{\nu}_\rho$ is at most

$$\rho^k((1+\rho+\rho^2+\cdots)+(\overline{q}+\overline{q}^2+\cdots)) = \rho^k\left(\frac{1}{1-\rho}+\frac{\overline{q}}{1-\overline{q}}\right),$$

where the first series bounds the probability that the discrepancy at $-n-k$ moves after the clock rings at some site $y \leq -n-k$, while the second series bounds the probability that the discrepancy at $-n-k$ moves after the clock rings at some site $y > -n-k$.

Negative discrepancies at sites $y > n$ may contribute because the clock rings at their site or because the clock rings at a site $z$ where there is an $\eta$ and a $\xi$ particle. In this last case both the $\eta$ and $\xi$ particles may reach $y$ forcing the negative discrepancy at $y$ to move. For this to happen all the sites from $z$ to either $y-1$ or $y+1$ have to be occupied by $\eta$ particles, which under $\widetilde{\nu}_\rho$ are distributed as a Bernoulli product measure of parameter $\rho$. Thus, the contribution of a negative discrepancy at $y = n+k$ to $\int \widetilde{L}^+ g_n \, d\widetilde{\nu}_\rho$ is at most

$$\overline{q}^k(1+2\rho+2\rho^2+\cdots) < 2\overline{q}^k\left(\frac{1}{1-\rho}\right).$$

We now choose $k_0$ large enough to satisfy

$$\sum_{k \geq k_0} \rho^k\left(\frac{1}{1-\rho}+\frac{\overline{q}}{1-\overline{q}}\right) + 2\overline{q}^k\left(\frac{1}{1-\rho}\right) < \frac{\delta_2}{2}.$$

Hence the positive terms due to discrepancies within distance $k_0$ of $[-n,n]$ contribute at least $\delta_2/2$. Note that if these terms are present for a given $(\eta,\xi)$ configuration, then $g_{n+k_0}(\eta,\xi) \geq g_n(\eta,\xi) + 1$. Let $s_x(\eta,\xi)$ be the length of the longest interval containing $x$ whose sites are all occupied by $\eta$ particles. Then, $s_x(\eta,\xi)$ is an upper bound for the rate at which a positive discrepancy at $x$ starts to move. Therefore the contribution of the positive discrepancies in the interval $[-n-k_0,-n-1]$ to $\int \widetilde{L}^+ g_n \, d\widetilde{\nu}_\rho$ is at most

$$\int \sum_{x=-n-k_0}^{-n-1} s_x(\eta,\xi) \mathbf{I}_{\{g_{n+k_0}(\eta,\xi) > g_n(\eta,\xi)\}} \, d\widetilde{\nu}_\rho.$$

Similarly, the contribution of negative discrepancies in the interval $[n+1, n+k_0]$ is bounded above by

$$\int \sum_{x=n+1}^{n+k_0} s'_x(\eta,\xi) \mathbf{I}_{\{g_{n+k_0}(\eta,\xi) > g_n(\eta,\xi)\}} \, d\widetilde{\nu}_\rho,$$

where $s'_x(\eta,\xi)$ is the length of the longest interval containing $x$ such that all its sites distinct from $x$ are occupied by $\eta$ particles. Hence, at least one of the two integrals above is bounded below by $\delta_2/4$. Since under $\widetilde{\nu}_\rho$ both $\sum_{x=-n-k_0}^{-n-1} s_x(\eta,\xi)$ and $\sum_{x=n+1}^{n+k_0} s'_x(\eta,\xi)$ are integrable and have distributions that do not depend on $n$, we conclude that there exists a $\delta_3 > 0$ such that for all $n \geq n_0$ we have

$$\widetilde{\nu}_\rho(g_{n+k_0} \geq g_n + 1) \geq \delta_3 > 0.$$

This implies that $\int g_n \, d\widetilde{\nu}_\rho$ grows linearly with $n$. It now follows from Lemmas 6.1 and 6.4 that $\int \widetilde{L}^- f_n \, d\widetilde{\nu}_\rho$ diverges as $n$ tends to infinity. Therefore, to get a contradiction it suffices to show that $\int \widetilde{L}^+ f_n \, d\widetilde{\nu}_\rho$ is bounded above by a constant that does not depend on $n$. This is done as follows: Arguing as we did for $\int \widetilde{L}^+ g_n \, d\widetilde{\nu}_\rho$ we see that a positive discrepancy at $-n-k$ contributes at most

$$\rho^k((1+\rho+\rho^2+\cdots)+(\overline{q}+\overline{q}^2+\cdots)) = \rho^k\left(\frac{1}{1-\rho}+\frac{\overline{q}}{1-\overline{q}}\right),$$

and that a positive discrepancy at $n+k$ contributes at most

$$\overline{q}^k(1+2\rho+2\rho^2+\cdots) < 2\overline{q}^k\left(\frac{1}{1-\rho}\right).$$

Therefore, we obtain as upper bound:

$$\sum_{k=1}^\infty \rho^k\left(\frac{1}{1-\rho}+\frac{\overline{q}}{1-\overline{q}}\right) + \sum_{k=1}^\infty 2\overline{q}^k\left(\frac{1}{1-\rho}\right),$$

which is finite and does not depend on $n$. □

PROOF OF THEOREM 1.6.   Let $\mu \in \mathcal{I}$ be such that $\mu(\{\mathbf{1}\}) = 0$ and let $\widetilde{\nu}_\rho$ be as in Proposition 6.1. Discrepancies move at rate at least 1 and the sites they visit are part of the sites visited by a $p,q$ random walk. Therefore, if there are only a finite number of positive (negative) discrepancies to the left of the origin, then the probability of having a positive (negative) discrepancy at a given site $x$ tends to 0 as time goes to infinity. To see this for a given configuration $(\eta,\xi)$, use an approximating sequence of finite configurations $(\eta_n,\xi_n)$ and note that the probability of having a positive (negative) discrepancy at a given site goes to 0 uniformly in $n$ as time goes to infinity. Hence, it follows from Proposition 6.1 and the invariance of $\widetilde{\nu}_\rho$ that for all $x < y$

$$\widetilde{\nu}_\rho\{(\eta,\xi):\eta(x)=\xi(y)\neq\eta(y)=\xi(x)\}=0,$$

therefore

$$\widetilde{\nu}_\rho\{(\eta,\xi):\eta\geq\xi \text{ or } \xi\geq\eta\}=1.$$

The theorem now follows from the remark at the end of Section 4.

The special case in which $p = 1$ needs an adaptation of the proof of Lemma 3.2 where we assumed that $p(x, y)$ was irreducible. In this particular case one can easily see that on the set $\int_0^1 L^{++} f_x(\eta_s^\xi)\,ds = \infty$ we also have $\int_0^1 L^{++} f_y(\eta_s^\xi)\,ds = \infty$ a.s. for all $y < x$. This observation allows us to prove the lemma with basically the same arguments as before. $\square$

**Acknowledgments.** The authors thank the referees for their comments and suggestions. H. Guiol thanks Pablo Ferrari and Tom Mountford for valuable discussions. H. Guiol thanks the Instituto de Matemática e Estatística of São Paulo University and the Institut de Mathématique of the Swiss Federal Institute of Technology at Lausanne for their wonderful hospitality during which part of this work has been achieved.

## REFERENCES

[1] ANDJEL, E. D. (1981). The asymmetric simple exclusion process on $Z^d$. *Z. Wahrsch. Verw. Gebiete* **58** 423–432. MR639150
[2] FELLER, W. (1971). *An Introduction to Probability Theory and Its Applications* **2**, 2nd ed. Wiley, New York. MR0270403
[3] GUIOL, H. (1997). Un résultat pour le processus d'exclusionà longue portée. [A result for the long-range exclusion process.] *Ann. Inst. H. Poincaré Probab. Statist.* **33** 387–405. MR1465795
[4] KESTEN, H. (1972). Sums of independent random variables—without moment conditions. *Ann. Math. Statist.* **43** 701–732. MR301786
[5] LIGGETT, T. M. (1976). Coupling the simple exclusion process. *Ann. Probab.* **4** 339–356. MR418291
[6] LIGGETT, T. M. (1980). Long range exclusion processes. *Ann. Probab.* **8** 861–889. MR586773
[7] LIGGETT, T. M. (1985). *Interacting Particle Systems*. Springer, New York. MR776231
[8] SPITZER, F. (1970). Interaction of Markov processes. *Adv. Math.* **5** 246–290. MR268959
[9] ZHENG, X. G. (1988). Ergodic theorem for generalized long-range exclusion processes with positive recurrence transition probabilities. *Acta Math. Sinica (N.S.)* **4** 193–209. MR965568

LATP-CMI
36 RUE JOLIOT-CURIE
13013 MARSEILLE
FRANCE
E-MAIL: enrique.andjel@cmi.univ-mrs.fr

TIMC-TIMB
PAVILLON D
FACULTÉ DE MÉDECINE
38706 LA TRONCHE–CEDEX
FRANCE
E-MAIL: herve.guiol@imag.fr